\documentclass{gtart_h}  

\def\ifplaintex{\expandafter\ifx\csname documentclass\endcsname\relax}

\def\ifplaintex{\expandafter\ifx\csname documentclass\endcsname\relax}

\ifplaintex 
\else
\fi

\expandafter\ifx\csname epsfbox\endcsname\relax\input epsf\fi

\def\gt{{\mathsurround=0pt\it $\cal G\mskip-2mu$eometry \&\ 
$\cal T\!\!$opology}}        

\def\gtp{{\mathsurround=0pt\it $\cal G\mskip-2mu$eometry \&\ 
$\cal T\!\!$opology $\cal P\!$ublications}}  

\def\lognumber#1{\def\thelognumber{#1}}
\def\volumenumber#1{\def\thevolumenumber{#1}}
\def\papernumber#1{\def\thepapernumber{#1}}
\def\volumeyear#1{\def\thevolumeyear{#1}}

\def\pagenumbers#1#2{\def\startpage{#1}\def\finishpage{#2}}
\def\published#1{\def\publishdate{#1}}
\def\proposed#1{\def\theproposer{#1}}
\def\seconded#1{\def\theseconders{#1}}
\def\received#1{\def\receiveddate{#1}}
\def\revised#1{\def\reviseddate{#1}}
\def\accepted#1{\def\accepteddate{#1}}
\def\asciititle#1{\def\theasciititle{#1}}
\def\covertitle#1{\def\thecovertitle{#1}}

\def\asciiaddress#1{\def\theasciiaddress{#1}}
\def\asciiemail#1{\def\theasciiemail{#1}}

\long\def\asciiabstract#1{\long\def\theasciiabstract{#1}}

\let\\\par\let\thelognumber\relax
\let\thevolumenumber\relax\let\thepapernumber\relax
\let\thevolumeyear\relax\let\thesamplenumber\relax\let\startpage\relax
\let\finishpage\relax\let\publishdate\relax\let\receiveddate\relax
\let\reviseddate\relax\let\accepteddate\relax\let\theasciititle\relax
\let\thecovertitle\relax\let\theasciiauthors\relax\let\theasciiaddress\relax
\let\theasciiabstract\relax
\let\theasciiemail\relax\let\theshortauthors\relax\let\theshorttitle\relax

\long\def\maketitlep{   

\count0=\startpage

\gt\hfill      
\hbox to 77pt{\vbox to 0pt{\vglue -15pt\epsfbox{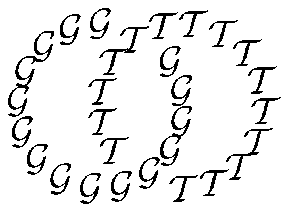}\vss}\hss}
\break
{\small\ifx\thesamplenumber\relax 
Volume \else Sample
\fi\thevolumenumber\ (\thevolumeyear)
\startpage--\finishpage\nl
Published: \publishdate}
\vglue 0.5truein plus 0.4fil minus 0.1truein

{\parskip=0pt\leftskip 0pt plus 1fil\def\\{\par\smallskip}{\ifplaintex\large
\else\Large\fi\bf\thetitle}\par\medskip}   

\vglue 0pt plus 0.1fil 

{\parskip=0pt\leftskip 0pt plus 1fil\def\\{\par}{\sc\theauthors}
\par\medskip}

\vglue 0pt plus 0.1fil 

{\small\parskip=0pt\let\newline\\
{\leftskip 0pt plus 1fil\def\\{\par}{\sl\theaddress}\par}
\expandafter\ifx\theemail\relax    
\relax\else\vglue 5pt plus 0.02fil minus 2pt\def\\{\stdspace{\rm 
and}\stdspace} 
\cl{Email:\stdspace\tt\theemail}\fi
\ifx\theurl\relax                  
\relax\else\vglue 5pt plus 0.02fil minus 2pt\def\\{\stdspace{\rm 
and}\stdspace}
\cl{URL:\stdspace\tt\theurl}\fi\par}

\vglue 7pt plus 0.3fil minus 3pt

{\bf Abstract}
\vglue 5pt plus 0.1fil minus 2pt

\theabstract

\vglue 7pt plus 0.3fil minus 3pt

{\bf AMS Classification numbers}\quad Primary:\quad \theprimaryclass

Secondary:\quad \thesecondaryclass

\vglue 5pt plus 0.3fil minus 2pt

{\bf Keywords:}\quad \thekeywords

\vglue 10pt plus 0.5fil minus 5pt

{\small  Proposed: \theproposer\hfill Received: \receiveddate\nl
Seconded: \theseconders\hfill 
\ifx\reviseddate\relax                         
Accepted: \accepteddate                        
\else
Revised: \reviseddate                          
\fi}
\eject
}       

\font\phead=cmsl9 scaled 950
\font=cmsl9 scaled 1050
\font\pnum=cmbx10 scaled 913
\font\lnum=cmbx10 
\font\pfoot=cmsl9 scaled 950
\font=cmsl9 scaled 1050
\ifplaintex
\headline{\vbox to 0pt{\vskip -4.5mm\line{\small\phead\ifnum
\count0=\startpage ISSN 1364-0380 (on line)
1465-3060 (printed) \hfill {\pnum\folio}\else\ifodd\count0\def\\{ }%
\ifx\theshorttitle\relax\thetitle\else\theshorttitle\fi\hfill{\pnum\folio}
\else\def\\{ and }{\pnum\folio}\hfill\ifx\theshortauthors\relax\theauthors
\else\theshortauthors\fi\fi\fi}\vss}}
\footline{\vbox to 0pt{\vglue 0mm\line{\small\pfoot\ifnum\count0=\startpage
\copyright\ \gtp\hfill\else
\gt, Volume \thevolumenumber\ (\thevolumeyear)\hfill\fi}\vss
}}
\else
\makeatletter
\def\@oddhead{{\small\lhead\ifnum\count0=\startpage ISSN 1364-0380 (on line)
1465-3060 (printed) \hfill {\lnum\number\count0}\else\ifodd\count0
\def\\{ }\ifx\theshorttitle\relax \thetitle \else\theshorttitle\fi\hfill
{\lnum\number\count0}\else\def\\{ and }{\lnum\number\count0}
\hfill\ifx\theshortauthors\relax 
\theauthors\else\theshortauthors\fi\fi\fi}}\def\@evenhead{\@oddhead}
\def\@oddfoot{\small\lfoot\ifnum\count0=\startpage\copyright\ \gtp\hfill\else
\gt, Volume \thevolumenumber\ (\thevolumeyear)\hfill\fi}
\def\@evenfoot{\@oddfoot}
\makeatother
\fi

\newwrite\gtoutfile
\long\gdef\makeheadfile{  
{\def\\{, }\def\s{ }
\immediate\openout\gtoutfile head.xxx
\immediate\write\gtoutfile{Proxy-for: \ifx\theasciiauthors\relax
\theauthors\else\theasciiauthors\fi\s<\ifx\theasciiemail\relax\theemail\else\theasciiemail\fi>}
\immediate\write\gtoutfile{\noexpand\\}
\immediate\write\gtoutfile{Authors: \ifx\theasciiauthors\relax
\theauthors\else\theasciiauthors\fi}
{\def\\{ }\immediate\write\gtoutfile{Title: \ifx\theasciititle\relax
\thetitle\else\theasciititle\fi}}
\immediate\write\gtoutfile{Subj-class: GT or SG or MG etc}
\immediate\write\gtoutfile{MSC-class: \theprimaryclass\ifx\thesecondaryclass\relax\else, \thesecondaryclass\fi}
\immediate\write\gtoutfile{Journal-ref: Geom. Topol. \thevolumenumber
(\thevolumeyear) \startpage-\finishpage}
\immediate\write\gtoutfile{Comments: Published by Geometry and Topology at}
\immediate\write\gtoutfile{\s\s http://www.maths.warwick.ac.uk/gt/GTVol\thevolumenumber/paper\thepapernumber.abs.html}
\immediate\write\gtoutfile{\noexpand\\}
\immediate\write\gtoutfile{}
\ifx\theasciiabstract\relax
\immediate\write\gtoutfile{\theabstract}\else
\immediate\write\gtoutfile{\theasciiabstract}\fi
\immediate\write\gtoutfile{}
\immediate\write\gtoutfile{\noexpand\\}
\immediate\write\gtoutfile{}
\immediate\closeout\gtoutfile}}  

\def\maketitlepage{\maketitlep\makeheadfile}
\let\maketitle\maketitlepage

\lognumber{407}
\received{15 January 2004}
\volumenumber{8}\papernumber{38}\volumeyear{2004}
\pagenumbers{1385}{1425}
\revised{1 September 2004}
\published{27 October 2004}
\accepted{11 October 2004}
\proposed{Bill Dwyer}
\seconded{Thomas Goodwillie, Gunnar Carlsson}

\usepackage{amsmath,amssymb,amscd}
\usepackage[all]{xy}
\usepackage[mathscr]{euscript}

\RequirePackage{xspace}

\usepackage{enumerate}        
\newenvironment{roenumerate}{\begin{enumerate}[\upshape (i)]}{\end{enumerate}}

\newcommand\nc {\newcommand}
\newcommand\rnc{\renewcommand}
\nc\script{\mathscr}

\newcount\blopone
\newcount\xone
\newcount\xtwo
\newcount\ytwo

\newtheorem{theorem}{Theorem}[section]
\newtheorem{prop}[theorem]{Proposition}
\newtheorem{equivform}[theorem]{Equivalent Formulation}

\newtheorem{refinement}[theorem]{Refinement}
\newtheorem{summary}[theorem]{Summary}
\newtheorem{importnota}[theorem]{Important Notation}
\newtheorem{prblm}[theorem]{Problem}
\newtheorem{notation}[theorem]{Notation}
\newtheorem{defin}[theorem]{Definition}
\newtheorem{caution}[theorem]{Caution}
\newtheorem{remark}[theorem]{Remark}
\newtheorem{reminder}[theorem]{Reminder}
\newtheorem{lemma}[theorem]{Lemma}
\newtheorem{construction}[theorem]{Construction}
\newtheorem{corollary}[theorem]{Corollary}
\newtheorem{example}[theorem]{Example}
\newtheorem{conclusion}[theorem]{Conclusion}
\newtheorem{triviality}[theorem]{Triviality}
\newtheorem{proto}[theorem]{Prototype Quasifibration}
\newtheorem{cauex}[theorem]{Cautionary Example}
\newtheorem{hypo}[theorem]{Hypothesis}
\newtheorem{subth}{ }[theorem]
\newtheorem{case}{Case}[theorem]
\newtheorem{ssubth}{ }[subth]

\nc\tri[1]{\begin{triviality} \label{#1}}
\nc\eqf[1]{\begin{equivform} \label{#1}}
\nc\cas[1]{\begin{case} \label{#1}}
\nc\rfn[1]{\begin{refinement} \label{#1}}
\nc\prt[1]{\begin{proto} \label{#1}}
\nc\lem[1]{\begin{lemma} \label{#1}}
\nc\pro[1]{\begin{prop} \label{#1}}
\nc\thm[1]{\begin{theorem} \label{#1}}
\nc\cor[1]{\begin{corollary} \label{#1}}
\nc\dfn[1]{\begin{defin}\label{#1}\begin{em}}
\nc\sthm[1]{\begin{subth} \label{#1}}
\nc\exm[1]{\begin{example} \label{#1}\begin{em}}
\nc\plm[1]{\begin{prblm} \label{#1}}
\nc\rmk[1]{\begin{remark} \label{#1}\begin{em}}
\nc\rmd[1]{\begin{reminder} \label{#1}\begin{em}}
\nc\ntn[1]{\begin{notation} \label{#1}\begin{em}}
\nc\smr[1]{\begin{summary} \label{#1}\begin{em}}
\nc\cau[1]{\begin{caution} \label{#1}}
\nc\hyp[1]{\begin{hypo} \label{#1}}
\nc\imn[1]{\begin{importnota} \label{#1}}
\nc\cax[1]{\begin{cauex} \label{#1}}
\nc\con[1]{\begin{construction} \label{#1}}
\nc\ssthm[1]{\begin{ssubth} \label{#1}}
\nc\cnc[1]{\begin{conclusion} \label{#1}}

\nc\elem{\end{lemma}}
\nc\eeqf{\end{equivform}}
\nc\erfn{\end{refinement}}
\nc\eprt{\end{proto}}
\nc\ethm{\end{theorem}}
\nc\ecor{\end{corollary}}
\nc\edfn{\end{em}\end{defin}}
\nc\esthm{\end{subth}}
\nc\epro{\end{prop}}
\nc\etri{\end{triviality}}
\nc\eexm{\end{em}\end{example}}
\nc\ermk{\end{em}\end{remark}}
\nc\ermd{\end{em}\end{reminder}}
\nc\eplm{\end{prblm}}
\nc\ecas{\end{case}}
\nc\ecau{\end{caution}}
\nc\ecax{\end{cauex}}
\nc\eimn{\end{importnota}}
\nc\entn{\end{em}\end{notation}}
\nc\econ{\end{construction}}
\nc\esmr{\end{em}\end{summary}}
\nc\ehyp{\end{hypo}}
\nc\ecnc{\end{conclusion}}
\nc\essthm{\end{ssubth}}

\nc\sst{\scriptstyle}
\newcommand{\comment}[1]{}
\newcommand{\ri}{\longrightarrow}

\newcommand{\zz}{{\mathbb Z}}

\nc\bR{{\mathbf R}}
\nc\bS{{\mathbf S}}
\nc\bSR{{\mathbf S}_{\mathbf R}^{}}
\nc\bT{{\mathbf T}}
\nc\bTd{{\mathbf T'}}
\nc\bU{{\mathbf U}}
\nc\z{\zeta}
\nc\bc{{\mathbb{BC}}}
\nc\ct{{\script T}}
\nc\cs{{\script S}}
\nc\car{{\script R}}
\nc\ca{{\script A}}
\nc\cb{{\script B}}
\nc\cc{{\script C}}
\nc\cd{{\script D}}
\nc\ce{{\script E}}
\nc\ci{{\script I}}
\nc\cu{{\script U}}
\nc\bZ{{\mathbb Z}}
\nc\bd{\begin{description}}
\nc\ed{\end{description}}
\nc\ctob{{\script C}at\big(\ci^{op},\ca\big)}
\nc\clim{{\ds\mathop{\rm lim}_{\ds\longleftarrow}}}
\nc\climi{\clim^{\!i}\,}
\nc\climn{\clim^{\!n}\,}
\nc\colim{{\ds\mathop{\rm colim}_{\ds\la}}}
\nc\oa{\overline{\ca}}
\nc\s{\sigma}
\nc\ta{\tau}
\nc\os{\overline\sigma}
\nc\ot{\overline\tau}
\nc\T{\Sigma}
\nc\de[1]{{\mathop{\rm deg(#1)}}}
\nc\Ad[1]{\mathop{\rm Ad}(#1)}
\nc\ad[1]{\mathop{\rm ad}(#1)}
\nc\Tor{\text{\rm Tor}}
\nc\be{\begin{roenumerate}}
\nc\ee{\end{roenumerate}}

\def\der #1 {D\left(#1\right)}

\nc\prf{\begin{proof}}
\nc\eprf{\end{proof}}
\nc\ds{\displaystyle}

\nc\ab{{\script A}b}

\nc\csab{{\script C}at\big(\cs^{op},\ab\big)}
\nc\ctab{{\script C}at\Big({\{\ct^\alpha\}}^{op},\ab\Big)}
\nc\csex{{\script E}x\big(\cs^{op},\ab\big)}
\nc\ctex{{\script E}x\Big({\{\ct^\alpha\}}^{op},\ab\Big)}
\nc\sub{\qquad\subset\qquad}
\nc\ctr[1]{{\left.\ct\left(-,#1\right)\right|}_{\cs}}
\nc\ctrf[2]{{\left.\ct\left(#1,#2\right)\right|}_{\cs}}
\nc\Ctr[1]{{\left.\ct\left(-,#1\right)\right|}_{\ct^\alpha}}
\nc\Ctrf[2]{{\left.\ct\left(#1,#2\right)\right|}_{\ct^\alpha}}

\nc\la{\longrightarrow}
\nc\oti{{^L\otimes_A^{}}}
\nc\rs{\s^{-1}A}
\nc\rrs{{\{\s^{-1}A\}}\op}
\nc\br{{\{\s^{-1}A\}}}

\nc\nin{\noindent}
\nc\cad[1]{\text{card}(#1)}
\nc\eq{\quad=\quad}
\nc\BA{\begin{array}{c}}
\nc\EA{\end{array}}
\nc\kth{$K$--theory}

\nc\barr{
\[
\begin{array}{cccccccccccccccc}
}
\nc\earr{
\end{array}
\]
}
\nc\as[1]{{\langle S\rangle}^{#1}}

\nc\yy[1]{{\left.\ct\left(-,#1\right)\right|}_{\ct^c}}
\nc\vrep[2]{{\left.\ct\left(#1,#2\right)\right|}_{\ct^\alpha}}
\nc\da{\downarrow}
\nc\Hom{{\mathop{\rm Hom}}}
\nc\End{{\mathop{\rm End}}}
\nc\Ext{{\mathop{\rm Ext}}}
\nc\mr\Modtc
\nc\PExt{{\mathop{\rm PExt}}}
\nc\bA{{\mathbf A}}
\nc\bB{{\mathbf B}}
\nc\bC{{\mathbf C}}
\nc\bD{{\mathbf D}}

\nc\ra\ri
\nc\Mod[1]{\ensuremath{\mathop{\textup{Mod-}#1}}\xspace}
\nc\Md {\ensuremath{\mathop{\textup{Mod}}}}
\rnc\mod[1]{\ensuremath{\mathop{\textup{mod-}#1}}\xspace}
\nc\Modtc{\Mod{\ct^c}}
\nc\pgldim[1]{\mathop{\rm pgldim}\,#1}
\nc\tstr{$t$--structure}
\nc\perf{^{\hbox{\rm\tiny perf}}}
\nc\op{^{\hbox{\rm\tiny op}}}
\nc\p{\varphi}

\begin{document}

\authors{Amnon Neeman\\Andrew Ranicki}
\address{Centre for Mathematics and its Applications, 
The Australian National University\\Canberra, ACT 0200, 
Australia\\\smallskip\\School of Mathematics, University of 
Edinburgh\\Edinburgh EH9 3JZ, Scotland, UK}

\gtemail{\mailto{Amnon.Neeman@anu.edu.au}, \mailto{a.ranicki@ed.ac.uk}}
\asciiemail{Amnon.Neeman@anu.edu.au, a.ranicki@ed.ac.uk}
\asciiaddress{Centre for Mathematics and its Applications, 
The Australian National University\\Canberra, ACT 0200, 
Australia\\and\\School of Mathematics, University of 
Edinburgh\\Edinburgh EH9 3JZ, Scotland, UK}

\title{Noncommutative localisation in algebraic $K$--theory I}
\covertitle{Noncommutative localisation\\in algebraic $K$--theory I}
\asciititle{Noncommutative localisation in algebraic K-theory I}

\begin{abstract}
This article establishes, for an appropriate localisation of associative
rings, a long exact sequence in algebraic \kth. The main result
goes as follows. Let $A$ be an associative ring and let 
$A \ri B$ be the
localisation with respect to a set $\sigma$ of maps between finitely
generated projective $A$--modules. Suppose that $\Tor_n^A(B,B)$ vanishes for
all $n>0$. View each map in $\sigma$ as a complex (of length 1, meaning 
one non-zero map between two non-zero objects) in the
category of perfect complexes $D\perf(A)$. Denote by $\langle\sigma\rangle$ the
thick subcategory generated by these complexes. 
Then the canonical
functor $D\perf(A)\ri D\perf(B)$ induces (up to direct factors) an
equivalence
${\ds D\perf(A)}/{\ds\langle\sigma\rangle} \ri D\perf(B)$.
As a consequence, one
obtains a homotopy fibre sequence 
$$\CD
K(A,\sigma) @>>> K(A)
@>>> K(B)
\endCD$$
(up to surjectivity of $K_0(A)\ri K_0(B)$) of Waldhausen \kth\ 
spectra.

In subsequent articles~\cite{Neeman-Ranicki03,Neeman-Ranicki04}
we will present the $K$-- and $L$--theoretic consequences
of the main theorem in a form more suitable for the applications to
surgery. For example if, in addition to the vanishing of $\Tor_n^A(B,B)$,
we also assume that every map in $\s$ is a monomorphism, then
there is a description of the homotopy fiber of the map $K(A)\ri K(B)$
as the Quillen \kth\ of a suitable exact category of torsion modules.
\end{abstract}

\asciiabstract{
This article establishes, for an appropriate localisation of
associative rings, a long exact sequence in algebraic K-theory. The
main result goes as follows. Let A be an associative ring and let
A-->B be the localisation with respect to a set sigma of maps between
finitely generated projective A-modules. Suppose that Tor_n^A(B,B)
vanishes for all n>0. View each map in sigma as a complex (of length
1, meaning one non-zero map between two non-zero objects) in the
category of perfect complexes D^perf(A). Denote by <sigma> the thick
subcategory generated by these complexes.  Then the canonical functor
D^perf(A)-->D^perf(B) induces (up to direct factors) an equivalence
D^perf(A)/<sigma>--> D^perf(B).  As a consequence, one obtains a
homotopy fibre sequence K(A,sigma)-->K(A)-->K(B) (up to surjectivity
of K_0(A)-->K_0(B)) of Waldhausen K-theory spectra.

In subsequent articles we will present the K- and L-theoretic
consequences of the main theorem in a form more suitable for the
applications to surgery. For example if, in addition to the vanishing
of Tor_n^A(B,B), we also assume that every map in sigma is a
monomorphism, then there is a description of the homotopy fiber of the
map K(A)-->K(B) as the Quillen K-theory of a suitable exact category
of torsion modules.}


{\small\maketitle}

\section*{Introduction}\addcontentsline{toc}{section}{Introduction}
\label{S0}

Probably the most useful technical tool in algebraic \kth\ is localisation.
Localisation tells us about 
certain long exact sequences in \kth.
Given a ring $A$, Quillen told us how to assign to it a \kth\ 
spectrum $K(A)$.
Given rings $A$ and $\rs$, where
$\rs$ is a localisation of $A$, there is a map of spectra
$K(A)\ri K(\rs)$. A localisation theorem expresses the homotopy
fiber as $K(\bR)$, the \kth\ of some suitable $\bR$. In the
early localisation theorems $K(\bR)$ was Quillen's \kth\ of
some exact category of torsion modules. But in more recent,
more general theorems one allows $K(\bR)$ to be the Waldhausen
\kth\ of some suitable Waldhausen model category $\bR$.
In more concrete terms we get a long exact sequence
\[
\cdots \ri K_1(A)  \ri K_1(\rs) \ri
K_0(\bR) \ri K_0(A)  \ri K_0(\rs) \ri0.
\]
There has been extensive literature over the years, proving localisation
theorems in algebraic \kth. Let us
provide a brief sample of the existing literature.
For commutative rings, or more generally for
schemes, the reader can see
Bass'~\cite{Bass},
Quillen's~\cite{Quillen1}, Grayson's~\cite{Grayson76}, 
Levine's~\cite{Levine88}, Weibel's~\cite{Weibel84,Weibel89} 
and Thomason's~\cite{ThomTro}. For non-commutative rings see 
Grayson's~\cite{Grayson80}, Schofield's~\cite{Schofield},
Weibel and Yao's~\cite{WeibelYao89}
and Yao's~\cite{Yao92}.

Some situations are very well understood. For example, let
$A$ be an associative ring with unit. 
Let $\sigma\subset A$ be a multiplicative set of elements in the center
$Z(A)$ of the ring $A$.  It is very classical to define $\rs$ as the
ring of fractions $a/s$, with $a\in A$ and $s\in\s$.  The ring $\rs$ is
called the {\em commutative localisation} of $A$ with respect to the
multiplicative set $\s$.  Note that the rings $A$ and $\rs$ are not
assumed commutative; the only commutativity hypothesis is that
$\s\subset Z(A)$.  If every element of $\s$ is a non-zero-divisor then
the existence of 
a localisation exact sequence is classical.

Over the years people have found more general
localisation theorems in the \kth\ 
of non-commutative rings. What we want to do in this 
article is treat localisation in the generality
that comes up in topology. In 
applications to higher
dimensional knot theory and surgery, the ring $A$ might
be the group ring of the fundamental group of some manifold, and the
set $\s$ almost never lies in the center of $A$. We first 
need to remind the reader of the generality in which 
the localisation of rings arises in topology.

Let us agree to some notation first. Our rings  are all 
associative and have units. 
Let $A$ be a ring. When we say ``$A$--module'' without an adjective,
we mean left $A$--module.

Let $A$ be any non-commutative ring. Let $\s$ be any set of maps
of finitely generated, projective $A$--modules. In symbols
\[
\s=\{s_i\co P_i\ri Q_i\mid\text{where }P_i,Q_i\text{ are f.g. projective}\}.
\]

\dfn{s-inve}
A ring homomorphism $A\ri B$ is called {\em $\s$--inverting} if, for
all \mbox{$s_i\co P_i\ri Q_i$} in $\s$, the map
\[
\CD
B\otimes_A^{}
P_i@>1\otimes_A^{} s_i>> B\otimes_A^{} Q_i
\endCD
\]
is an isomorphism.
\edfn

The collection of $\s$--inverting homomorphisms $A\ri B$ is naturally
a category. A morphism in this category is a commutative triangle
of ring homomorphisms:
\[
\xymatrix@R-10pt@C+5pt{
  & & B\ar[dd] \\
A \ar[rru]\ar[rrd] & & \\
 & & B'
}
\]
There is an old observation, due to 
Cohn~\cite{Cohn71} and Schofield~\cite{Schofield}, which 
says that the category of $\s$--inverting homomorphisms $A\ri B$
has an initial object. 

\dfn{Dcohn}
The initial object in the category of $\s$--inverting homomorphisms
is called the {\em Cohn localisation} or the
{\em universal localisation} of $A$ with respect to $\s$. In this
article it will be denoted $A\ri\rs$.
\edfn

See Vogel~\cite{Vogel1980,Vogel1982}, Farber and
Vogel~\cite{FarberVogel}, Farber and Ranicki~\cite{FarberRanicki},
and Ranicki~\cite{Ranicki1998,Ranicki2003}
for some of the applications of the algebraic $K$-- and
$L$--theory of Cohn localisations in topology.

In this article, a {\em perfect complex} of $A$--modules is a 
bounded complex of finitely generated projective modules. Let $C\perf(A)$
be the Waldhausen category of all perfect complexes. That is, the objects 
are the perfect complexes, the morphisms are the chain maps, the weak 
equivalences are the homotopy equivalences, and the cofibrations are the
degreewise split monomorphisms. Let 
$D\perf(A)$ be the associated homotopy
 category; the objects are still the perfect complexes, but the 
morphisms are homotopy equivalence classes of chain maps.

\rmk{strictperfect}
Our convention is slightly different from the standard one.
In the literature there is a distinction made between {\em perfect complexes}
and {\em strictly perfect complexes.} What we call a perfect complex is what,
elsewhere in the literature, is often referred to as a 
strictly perfect complex. We have almost no use for perfect complexes which
are not strictly perfect. For this reason we let the adverb ``strictly''
be understood.
\ermk

The set $\s$ of maps $s_i\co P_i\ri Q_i$ can be thought of as a set
of objects in either $C\perf(A)$ or $D\perf(A)$
(the two categories share the same set of objects). 
We simply take the chain complexes
\[
\CD
\cdots @>>> 0 @>>> P_i @>s_i>> Q_i @>>> 0 @>>> \cdots
\endCD
\]
We define a {\em Waldhausen category} $\bR\subset C\perf(A)$ as follows:

\dfn{Dwaldlift}
The category $\bR$ is the smallest subcategory of $C\perf(A)$
which
\be
\item Contains all the complexes in $\s$, in the sense above.
\item Contains all acyclic complexes.
\item Is closed under the formation of mapping cones and suspensions.
\item Contains any direct summand of any of its objects.
\ee
\edfn

The main {\it K}--theoretic result of the article becomes:

\thm{Tmainform}
Suppose $A$ is a ring, $\s$ a set of maps of finitely generated,
projective $A$--modules. 
Suppose $\Tor_n^A(\rs,\rs)=0$ for all $n>0$. 
Then the homotopy fiber of the map
$K(A)\ri K(\rs)$ is naturally identified, up to the failure 
of surjectivity of the map $K_0^{}(A)\ri K_0^{}(\rs)$,
with the spectrum $K(\bR)$. By $K(\bR)$ we mean the Waldhausen
\kth\ of the Waldhausen category $\bR$
of Definition~\ref{Dwaldlift}.
\ethm

Theorem~\ref{Tmainform} gives a precise version of what was stated,
slightly less precisely, in the abstract. It turns out that
Theorem~\ref{Tmainform} 
is a consequence of a statement about triangulated
categories. Next we explain the
triangulated categories results, and why the {\it K}--theoretic
statement in Theorem~\ref{Tmainform} is a formal 
consequence.  

Let $\car^c$ be the smallest triangulated subcategory of $D\perf(A)$,
containing $\s$ and containing all direct summands of any of its 
objects. The category $\bR$ of Definition~\ref{Dwaldlift} is
simply a Waldhausen model for $\car^c$.

\rmk{abs}
In the abstract we used the notation $\langle\s\rangle$ to denote
what we now call $\car^c$. The new notation is to avoid ambiguity,
clearly distinguishing the subcategory $\bR$
 generated by $\s$ in $C\perf(A)$ from the subcategory $\car^c$
generated by $\s$ in $D\perf(A)$.
\ermk

Let $\ct^c$ be the idempotent
completion of the Verdier quotient $D\perf(A)/\car^c$.
That is, form the Verdier quotient, and in it split all
idempotents. The main theorem, in its triangulated category
incarnation, asserts:

\thm{relat1}
Consider the natural functor $D\perf(A)\ri D\perf(\rs)$, which 
takes a complex in $D\perf(A)$ and tensors it with
$\rs$. Take the canonical factorisation
\[
\CD
D\perf(A)  @>\pi>>\ct^c @>T>>D\perf(\rs).
\endCD
\]
Then the following are equivalent:
\be
\item The functor $T\co\ct^c\ri D\perf(\rs)$ is an equivalence of categories.
\item For all $n\geq 1$ the group $\Tor_n^A(\rs,\rs)=0$.\footnote{The
    implication (i)$\Longrightarrow$(ii) may be 
found in Geigle and Lenzing~\cite{GeigleLenzing91}. The implication
    (ii)$\Longrightarrow$(i) is related to the telescope conjecture.
  The relation is studied in Krause~\cite{Krause03}.
    Unfortunately the telescope conjecture is false in general; see Keller~%
\cite{Keller94a}. More precisely: 
there are hardly any non-commutative rings for
which the telescope conjecture is known to hold. What we give here is a
    proof independent of the telescope conjecture.  
}
\ee
We call the localisation $A\ri\rs$ {\em stably flat}\footnote{In
  Dicks' article~\cite[p.\ 565]{Dicks77} the terminology for such a map
  is {\em a lifting,} while in
  Geigle and Lenzing's~\cite{GeigleLenzing91} it would be called a
  {\em homological epimorphism.} When we coined the term {\em stably
    flat} we were unaware of the earlier literature.} 
if the equivalent 
conditions above hold.
\ethm

Next we will sketch how Theorem~\ref{Tmainform}
follows from Theorem~\ref{relat1} (more detail will be given in later 
sections).
We have a sequence
\[
\CD
\car^c @>>>D\perf(A) @>>> 
D\perf(\rs).
\endCD
\]
This sequence always has a lifting to Waldhausen models
\[
\CD
\bR @>>>C\perf(A) @>>> C\perf(\rs).
\endCD
\]
Assume further that $\Tor_n^A(\rs,\rs)=0$ for all $n\geq1$.
Theorem~\ref{relat1} tells us that the category 
$D\perf(\rs)$ is 
identified, up to splitting idempotents, as the 
Verdier quotient ${D\perf(A)}/{\car^c}$.
Waldhausen's Localisation Theorem
(Theorem~\ref{Local}), coupled with Grayson cofinality
(Theorem~\ref{Graycof}), tell us
that
\[
\CD
K(\bR) @>>>K\big(C\perf(A)\big) @>\phi>> K\big(C\perf(\rs)\big)
\endCD
\]
identifies $K(\bR)$ as the $(-1)$--connected cover of the homotopy
fiber of the map $\phi$ above. By a theorem
of Gillet, the Waldhausen \kth\ $K(C\perf(A))$ is naturally isomorphic
to Quillen's $K(A)$; in other words, we get a commutative square where
the vertical maps are homotopy equivalences
\[
\CD
K(A) @>>> K(\rs) \\
@V|V\wr V @V|V\wr V \\
K\big(C\perf(A\big)) @>\phi>> K\big(C\perf(\rs)\big).
\endCD
\]
Hence $K(\bR)$ may be identified with the
 $(-1)$--connected cover of the homotopy
fiber of the map $K(A)\ri K(\rs)$. In other words: Theorem~\ref{Tmainform}
follows easily from Theorem~\ref{relat1}, modulo well known results of
Waldhausen, Grayson and Gillet.

At the level of homotopy groups this means we have an
infinite long exact sequence
$$\cdots\to K_1(\bR) \to K_1(A)  \to K_1(\rs)
  \to K_0(\bR) \to K_0(A)  \to K_0(\rs)$$
which continues indefinitely to the left. We are, however, {\em not} 
asserting that the map $K_0(A)  \ri K_0(\rs)$ is surjective.
The homotopy fiber $F$ of the map $K(A)\ri K(\rs)$ has 
in general a non-vanishing $\pi_{-1}^{}$, and the map $K(\bR)\ri F$
is an isomorphism in $\pi_i^{}$ for all $i\geq0$, but 
$\pi_{-1}^{}K(\bR)$
vanishes.

We have stated the main theorems mostly in the case
where the localisation is stably flat (see Theorem~\ref{relat1}
for the definition of stable flatness). 
There are examples of Cohn localisations which are not stably flat;
see~\cite{Neeman-Ranicki-Schofield}. Even in the non-stably-flat
 case the study of
the functor $T\co\ct^c\ri D\perf(\rs)$ is illuminating and has {\it
  K}--theoretic consequences. This article is devoted to studying the
functor $T$.

This article contains the proof of the theorems above, and
some
other formal, triangulated category facts about the functor
$T$. The applications
will appear separately. See~\cite{Neeman-Ranicki03} for {\it
  K}--theoretic
consequences, \cite{Neeman-Ranicki04} for consequences in {\it
  L}--theory,
as well as Krause's beautiful article~\cite{Krause03} which further
develops some of our results.

In presenting the proofs we tried to keep in mind that the reader might
not be an expert in derived categories. This is a paper of interest in
topology and surgery theory. We therefore try to give a
survey of the results, from the
literature on \kth\ and on triangulated categories,
which we need to appeal to. We give clear statements and 
careful references. We also try to
break down the proofs into a series of very easy steps.

The result is that the paper is much longer than necessary to
communicate the results to the experts; we ask the experts for
patience. The other drawback,
of presenting the proof in many easy steps, is that
the key issues can
become disguised. We will address this soon, 
in the discussion of the proof.

It is never clear how much the introduction ought to say about the
details of the proofs. Let us confine ourselves to the following. 
It is easy to produce the functors
\[
\CD
D\perf(A) @>\pi>> \ct^c @>T>> D\perf(\rs).
\endCD
\]
For any integer $n\in\zz$, they induce maps of abelian groups
\[
\CD
\Hom_{D\perf(A)}^{}(A,\T^n A) @. \\
@VVV @. \\
\Hom_{\ct^c}^{}(\pi A,\T^n \pi A)
@>{\p_n^{}}>> 
\Hom_{D\perf(\rs)}^{}(T\pi A,\T^n T\pi A).
\endCD
\]
If $T$ is an equivalence of categories, then the map $\p_n^{}$ above
must be an isomorphism. A minor variant of a theorem of Rickard's tells 
us that the converse also holds. To prove that $T$ is an equivalence
of categories, it suffices to show that the map $\p_n^{}$ is an isomorphism
for every $n\in\zz$. 

The construction of $T$ gives that $T\pi A=\rs$.
This means that the abelian group
$\Hom_{D\perf(\rs)}^{}(T\pi A,\T^n T\pi A)$ is just
$\Ext^n_{\rs}(\rs,\rs)$. It vanishes when $n\neq 0$.
For $n=0$, the endomorphisms of $\rs$, viewed
as a left $\rs$--module, are right multiplication by elements of
$\rs$. Therefore we are reduced to showing 
\[
\Hom_{\ct^c}^{}(\pi A,\T^n \pi A)=
\left\{
\begin{array}{lll}
0 &\qquad & \text{if }n\neq0\\
\rrs&\qquad & \text{if }n=0.
\end{array}
\right.
\]
In other words, the proof reduces to computing the groups
$\Hom_{\ct^c}^{}(\pi A,\T^n \pi A)$.

It happens to be very useful to turn the problem into one about 
unbounded complexes. Although Theorem~\ref{relat1} deals only with perfect 
complexes, the proof looks at $D(A)$, the unbounded derived category.
It is possible to embed the category $\ct^c$ in a larger category
$\ct$, and extend the map $\pi\co D\perf(A)\ri\ct^c$ to a map
$\pi\co D(A)\ri\ct$. What makes this useful
is that the extended functor $\pi$ has a right adjoint
$G\co \ct\ri D(A)$. By adjunction
\[
\Hom_{\ct^c}^{}(\pi A,\T^n \pi A)\eq
\Hom_{D(A)}^{}(A,\T^n G\pi A)\eq H^n(G\pi A),
\]
and we are reduced to computing $H^n(G\pi A)$.

\rmk{No torsion hyp}
It turns out that, for $n\geq0$, there is no need to assume
the vanishing of $\Tor_n^A(\rs,\rs)$. Without any hypotheses
we get $H^n(G\pi A)=0$ if $n>0$, while 
\[
H^0(G\pi A)\eq \Hom_{\ct^c}^{}(\pi A,\pi A)\eq \rrs.
\]
\ermk

The key lemma, which underpins everything we prove, is
Lemma~\ref{Gtrunc}. The lemma looks like a trivial little fact. It
asserts that, for the standard \tstr, the truncations of any object of
the form $G\pi x$ are also of the form $G\pi y$. This is the one
point where we use the fact that we are dealing with a Cohn
localisation, not just a general localisation in a triangulated
category. The lemma crucially depends on the complexes 
$\s$, which generate
the subcategory $\car^c$, being of length $\leq1$. 
In the case where $A$ is a commutative noetherian ring,
\cite{Neeman92B} tells us all the localisations of the 
derived category. 
It is easy to see that, without the hypothesis that
the complexes generating $\car^c$
be of length $\leq1$, essentially all our theorems fail.
The proof amounts to following the consequences of 
Lemma~\ref{Gtrunc}. We play around with some spectral sequences when
necessary, the argument is a little tricky at points, but none of this
changes the fact that Lemma~\ref{Gtrunc} is the foundation for
everything we prove.

This is our second attempt to expose the results; the first may be found in
\cite{Neeman-Ranicki01}. All but the experts in triangulated
categories found the first exposition difficult to read. As we have
already explained, this is our attempt to make the article readable.
We begin with a
survey of the main results we need from the literature. Then follows 
a sequence of easy steps, reducing the proof of Theorem~\ref{relat1}
to the computation of $H^n(G\pi A)$. The computations, which 
are the hard core of the article, come only at the end,
in sections~\ref{Sproofs},
\ref{n=0} and \ref{n<0}.

Since we want this article to be easy to read, we try not to assume that
the reader is very familiar with triangulated categories. We have therefore
gone to some trouble to keep our references to the literature focused.
In order to read the article, the triangulated category
background that is needed is:
\be
\item The standard \tstr\ on $D(A)$; \cite[Chapter 1]{BeiBerDel82}.
\item Homotopy limits and colimits; \cite[Sections 1--3]{Bokstedt-Neeman93}.
\item The generalisation of Thomason's localisation theorem;
\cite[Sections 1,2]{Neeman92A}
\ee
The reader will note that all the needed information is contained near
the beginning of the papers cited. We make a serious effort not to
refer any place else. But we feel free to quote any of the results
in the brief literature given above.

The fact that we cite only the three papers above leads to
historical inaccuracies. For example, 
the existence of the right adjoint $G$ 
to the functor $\pi\co D(A)\ri \ct$ was first proved by 
Bousfield~\cite{Bousfield75,Bousfield79A}. The many people who have
done excellent work in triangulated categories do not receive the
credit they deserve: see for example Keller's 
articles~\cite{Keller94,Keller94a} or Krause's~\cite{Krause98}.
Also,
there is a sense in which our main theorems are descended 
 from Thomason's~\cite{ThomTro}. In the survey article~\cite{Neeman05?}
we try to correct at least one of the historical inaccuracies,
indicating the crucial role of Thomason's work.

To keep the length from mushrooming to
infinity we have separated off the applications, which now appear in
\cite{Neeman-Ranicki03,Neeman-Ranicki04}.

{\bf Acknowledgements}\qua
The authors would like to thank the referees and editors
of \gt\ for many helpful comments, improving the
presentation of the results. The first author was partly supported 
by Australian Research Council grant DP0343239.

\section{Notation, and a reminder of {\it t}--structures}
\label{Ststr}

All our rings in this article will be associative rings with units.
Let $A$ be a ring.
Unless otherwise specified, all modules are left $A$--modules.
The derived category $D(A)$ means the unbounded
derived category of all complexes of $A$--modules. An object $x$
is a complex
\[
\CD
\cdots @>>> x^{n-2}_{} @>>>  x^{n-1}_{} 
@>\partial^{n-1}_{}>>  x^{n}_{} @>\partial^{n}_{}>>  x^{n+1}_{} @>>> 
 x^{n+2}_{} @>>> \cdots
\endCD
\]
As the reader has undoubtedly noticed, we write our complexes
cohomologically. Since one gets tired of adding a ``co'' to every word,
let it be understood. What we call
chain maps is what in the literature is
usually called cochain maps. What we call chain complexes 
is usually called cochain complexes.

The $n$th
homology of the complex $x$ above (which is what most people
refer to as the cohomology of the cochain complex) will be denoted
$H^n(X)$.

When it is clear which category we are dealing with, we write the
Hom--sets as $\Hom(x,y)$. When there are several categories around,
we freely use the notation $\ct(x,y)$ for $\Hom_\ct^{}(x,y)$.

In this article, {\it t}--structures on triangulated categories play a 
key role in many of our proofs. For an excellent exposition of
 this topic see Chapter 1 of \cite{BeiBerDel82}.  We give
here the bare essentials.
Let $\cs=D(A)$ as above. The only \tstr\ we use
in this article is the standard one on $\cs$. We remind the reader.

For any integer $n\in \zz$ there are two full subcategories of 
$\cs$. The objects are given by
\begin{eqnarray*}
\text{\rm Ob}(\cs^{\leq n}) &=& \{X\in\text{\rm Ob}(\cs)\mid
H^r(X)=0\,\,\hbox{\rm for all}~r>n\}, \\
\text{\rm Ob}(\cs^{\geq n})&=&\{X\in\text{\rm Ob}(\cs)\mid
H^r(X)=0\,\,\hbox{\rm for all}~r<n\}.
\end{eqnarray*}
The properties they satisfy are
\be
\item $\cs^{\leq n}\subset\cs^{\leq n+1}$.
\item $\cs^{\geq n}\subset\cs^{\geq n-1}$.
\item $\T\cs^{\geq n}=\cs^{\geq n-1}$, and $\T\cs^{\leq n}=\cs^{\leq n-1}$.
\item If $x\in \cs^{\leq -1}$ and $y\in \cs^{\geq 0}$, then
$
\Hom(x,y)=0.
$
\item For every object $x\in\cs$ there is a unique,
canonical distinguished triangle
\[
\CD
x^{\leq n-1} @>>> x @>>> x^{\geq n} @>>> \T x^{\leq n-1}
\endCD
\]
with $x^{\leq n-1}\in\cs^{\leq n-1}$ and $x^{\geq n}\in\cs^{\geq n}$.
\ee

\rmk{tsrtimp}
If $x$ is the complex
\[
\CD
\cdots @>>> x^{n-2}_{} @>>>  x^{n-1}_{} 
@>\partial^{n-1}_{}>>  x^{n}_{} @>\partial^{n}_{}>>  x^{n+1}_{} @>>> 
 x^{n+2}_{} @>>> \cdots
\endCD
\]
then the maps
\[
\CD
x^{\leq n-1} @>>> x @>>> x^{\geq n}
\endCD
\]
are concretely given by the chain maps
\[
\CD
\cdots \ri@.\quad @. x^{n-2}_{} 
@>>>  \hbox{Ker}(\partial^{n-1}_{}) @>>>  0 @>>> 0 @.\quad @.\ri
 \cdots \\
@. @. @VVV @VVV @VVV @VVV  @. @. \\
\cdots \ri@. @. x^{n-2}_{} @>>>  x^{n-1}_{} 
@>\partial^{n-1}_{}>>  x^{n}_{} @>>>  x^{n+1}_{} @. \quad @.\ri
 \cdots\\
@. @. @VVV @VVV @VVV @VVV  @. @. \\
\cdots\ri @. @. 0  @>>> 0  
@>>>  \hbox{Coker}(\partial^{n-1}_{}) @>>>   x^{n+1}_{} 
@. @. \ri 
 \cdots
\endCD
\]
\ermk

\section{Preliminaries, based on Waldhausen's work}
\label{Swaldh}
We begin with a brief review of
Waldhausen's foundational work. The
reader can find much more thorough treatments in
Waldhausen's article
\cite{Waldhausen85}, or in Section~1
of Thomason's~\cite{ThomTro}.

Let $\bS$ be a small
category with cofibrations and weak equivalences.
Out of $\bS$ Waldhausen constructs a spectrum, denoted 
$K(\bS)$. In Thomason's \cite{ThomTro} the category $\bS$
is assumed to be a full subcategory of the category of chain
complexes over some abelian category, the cofibrations
are maps of complexes which are split monomorphisms in each 
degree, and the weak equivalences contain the quasi-isomorphisms.
We will call such categories {\em permissible Waldhausen categories.}
In this article, we may assume that all categories with cofibrations
and weak equivalences are permissible Waldhausen categories.

\rmk{BiWald}
Thomason's term for them is {\em complicial biWaldhausen
categories.}
\ermk

\nin
Given a small, permissible Waldhausen category 
$\bS$, one can form its derived
category; just invert the weak equivalences. We denote this
derived category by $D(\bS)$. We have two major theorems
here, both of which are special cases of more general theorems
of Waldhausen. The first theorem may be found in 
Thomason's~\cite[Theorem 1.9.8]{ThomTro}:

\begin{theorem}[Waldhausen's Approximation Theorem] \label{Approx}
Let \mbox{$F\co\bS\ri\bT$} be an exact functor of small, permissible
Waldhausen categories (categories of chain complexes, as
above). Suppose that the induced map of derived categories
\[
D(F)\co D(\bS)\ri D(\bT)
\]
is an equivalence of categories. Then the induced map of
spectra
\[
K(F)\co K(\bS)\ri K(\bT)
\]
is a homotopy equivalence.
\end{theorem}

In this sense, Waldhausen's $K$--theory is almost an invariant
of the derived categories. To construct it one needs to
have a great deal more structure. One must begin with
a permissible category with cofibrations and weak equivalences.
But the Approximation Theorem asserts that the dependence
on the added structure is not strong.

Next we state Waldhausen's Localisation Theorem. 
The statement we give is an easy
consequence of
Theorem~\ref{Approx}, coupled with 
Waldhausen's~\cite[1.6.4]{Waldhausen85} or
Thomason's~\cite[1.8.2]{ThomTro}:

\begin{theorem}[Waldhausen's Localisation Theorem] \label{Local}
Let $\bR$, $\bS$ and $\bT$
be small, permissible
Waldhausen categories. Suppose 
\[
{\bR}\ri {\bS}\ri{\bT}
\] 
are exact functors of permissible
Waldhausen categories. Suppose further that 
\be
\item The induced
triangulated functors of derived categories
$$D({\bR})\ri D({\bS})\ri D({\bT})$$
compose to zero.
\item The functor $\p\co D({\bR})\ri D({\bS})$
is fully faithful.
\item If $x$ and $x'$ are objects of $D(\bS)$,
and the direct sum $x\oplus x'$ is
isomorphic in $D(\bS)$ to $\p(z)$ for some 
$z\in D(\bR)$, then $x,x'$ are isomorphic to $\p(y),\p(y')$ for
some $y,y'\in D(\bR)$.
\item
The natural map
\[
\CD
{\displaystyle D({\bS})}/{\displaystyle D({\bR})} @>>>
D({\bT})
\endCD
\] 
is an equivalence of categories. 
\ee
Then the sequence of spectra
 \[
K({\bR})\ri K({\bS})\ri K({\bT})
\]
is a homotopy fibration.
\ethm

We need one more general theorem, this one due to 
Grayson~\cite{Grayson87}.

\begin{theorem}[Grayson's Cofinality Theorem]\label{Graycof}
Let $\p\co\bT\ri\bTd$ be an exact functor of permissible
Waldhausen categories. Suppose the induced map $D(\p)\co D(\bT)\ri D(\bTd)$
is an idempotent completion. That is, the functor $D(\p)$ is fully faithful,
and for every object $x\in D(\bTd)$ there exists an object 
 $x'\in D(\bTd)$ and an isomorphism $x\oplus x'\simeq D(\p)(y)$, with
$y\in D(\bT)$.

Then the map of spectra
$K(\p)\co K(\bT)\ri K(\bTd)$ satisfies
\be
\item $K_i(\bT)\ri K_i(\bTd)$ is an isomorphism if $i\geq1$.
\item $K_0^{}(\bT)\ri K_0^{}(\bTd)$ is injective.
\ee
\ethm
In the article we will apply the
results of this section.
None of the results is very sensitive to changes in
Waldhausen models.
The additivity theorem, which we did not
discuss in this section, 
is sensitive to changes of permissible
Waldhausen categories. In this
article, and the two subsequent 
ones~\cite{Neeman-Ranicki03,Neeman-Ranicki04}, we never once use the 
additivity theorem. We can afford to confine ourselves
to proving the existence of one way to make the
choice of models.
Of course it is possible that, in the future, someone
will want to apply the results of the articles in conjunction with
the additivity theorem. Such a person will have to pay more
attention to the choice of Waldhausen categories.

Let us discuss one cheap way to produce models.

\lem{cheapmodels}
Let $\cs^c$ be a small triangulated category, $\car^c\subset\cs^c$
a triangulated subcategory containing all direct summands of
its objects. Suppose we are given a permissible
Waldhausen category $\bS$ and an equivalence of triangulated
categories
$\p\co D(\bS)\ri\cs^c$. Define $\bR$ to be the full Waldhausen
subcategory of all objects $x\in \bS$ so that $\p(x)$ 
is isomorphic in $\cs^c$ 
to an object in $\car^c\subset\cs^c$. Define the permissible
Waldhausen category $\bSR$ so that the objects, morphisms and
cofibrations are as in $\bS$, but the weak equivalences in $\bSR$
are the maps in $\bS$ whose mapping cones lie in $\bR$.

Then there is a commutative diagram of triangulated
functors, where the vertical maps are equivalences
\[
\CD
D(\bR) @>>>  D(\bS) @>>> D(\bSR) @.\\
@V|V\wr V       @V|V\wr V     @V|V\wr V @. \\
\car^c @>>> \cs^c  @>>> \cs^c/\car^c@..
\endCD
\]
\elem 

{\bf Idea of the Proof}\qua
The axioms of permissible
Waldhausen categories guarantee that the calculus of fractions,
in the passage from a permissible
Waldhausen category $\bS$ to its derived category $D(\bS)$, is
quite simple. Every morphism $x\ri y$ in $D(\bS)$
can be written as $\beta\alpha^{-1}_{}$, for
maps in $\bS$
\[
\CD
x @<\alpha<< x' @>\beta>> y
\endCD
\]
with $\alpha$ a weak equivalence. If $x$ and $y$ are in $\bR$ then,
since $\bR$ contains all the isomorphs in $D(\bS)$ of any of its objects,
$\bR$ must contain $x'$. Hence any morphism in $D(\bS)$ between
objects in the image of $D(\bR)$ lifts to $D(\bR)$. The equivalence
relation between pairs $(\alpha,\beta)$ as above is slightly more complicated,
but also only involves objects isomorphic in $\cs^c$ to $x$. Hence
$\beta\alpha^{-1}_{}$ will equal 
$\beta'{\{\alpha'\}}^{-1}_{}$
in $D(\bS)$ if and only if they are equal in $D(\bR)$. Thus the functor
$D(\bR)\ri D(\bS)$ is fully faithful. The objects in
$D(\bR)$ are, by definition of $\bR$, precisely the ones isomorphic 
to objects in $\car^c$.

The fact that $D(\bSR)\simeq\cs^c/\car^c$ is obvious from Verdier's 
construction of the quotient $\cs^c/\car^c$.
\hfill{$\Box$}

\section{The machine to produce examples}
\label{Sgenthom}

In order to apply the theorems of the last section,
we will produce
triangulated categories $\car^c\subset\cs^c$ and a 
triangulated functor $\cs^c/\car^c\ri\ct^c$ which is an idempotent
completion (as in Theorem~\ref{Graycof}).
There is a general machine which constructs examples. It is based
on a theorem by the first author. In this section
we will set up the notation, state the theorem and explain 
how it is applied.

\dfn{Dcompactobject}
Let $\cs$ be a triangulated category, containing all small coproducts
of its objects. An object $c\in\cs$ is called {\em compact} if every
map from $c$ to any coproduct factors through a finite part of the
coproduct. That is, any map
\[
\CD
c @>>> \ds\coprod_{\lambda\in\Lambda}t^{}_\lambda
\endCD
\]
factors as
\[
\CD
c @>>> 
\ds\coprod_{i=1}^nt^{}_{\lambda_i}
@>>> \ds\coprod_{\lambda\in\Lambda}t^{}_\lambda.
\endCD
\]
Equivalently, $c$ is compact if and only if
\[
\bigoplus_{\lambda\in\Lambda}\Hom(c,t^{}_\lambda)
\eq 
\Hom\left(c\,,\,\ds\coprod_{\lambda\in\Lambda}t^{}_\lambda\right).
\]
\edfn

\exm{Ecompact}
Let $A$ be a ring, and let $\cs=D(A)$ be the unbounded derived category
of $A$. The
category $\cs$ contains all small coproducts of its objects; we
can form direct sums of unbounded complexes.
Let $A\in D(A)=\cs$ be the chain complex
which is $A$ in degree $0$, and vanishes in all other degrees. 
For any $X\in\cs$ he have $\Hom(A,X)=H^0(X)$, and hence
$$\ds\bigoplus_{\lambda\in\Lambda}\Hom(A,t^{}_\lambda) =
\ds\bigoplus_{\lambda\in\Lambda}H^0(t^{}_\lambda)
  = H^0\left(\ds\coprod_{\lambda\in\Lambda}t^{}_\lambda\right)
  = \Hom\left(A\,,\,\ds\coprod_{\lambda\in\Lambda}t^{}_\lambda\right).$$
Thus $A$ is a compact object of $\cs$. 
\eexm

\dfn{Dcompactsub}
The full subcategory $\cs^c\subset\cs$ has for its objects all the
compact objects of $\cs$.
\edfn

\rmk{Rcomp}
It is easy to show that the full subcategory $\cs^c\subset\cs$ is closed
under triangles and direct summands.
\ermk

\exm{Ecomp}
In the special case of the category $\cs=D(A)$ of Example~\ref{Ecompact},
we know that $A$ is compact.
Any finite direct sum of
compact objects is compact, and any direct summand of a compact object 
is compact. We conclude that all finitely generated, 
projective $A$--modules are compact. The subcategory $\cs^c\subset\cs$
is triangulated, and hence we conclude that all bounded
chain complexes of finitely generated, projective modules are 
compact. In Corollary~\ref{Cmain.0.6} we will see that every
compact object in $\cs$ is isomorphic to a bounded complex 
of finitely generated, projective modules.
\eexm

\ntn{Notmainthm}
Next we set up the notation for the main theorem. 
Let $\cs$ be a triangulated category containing all small coproducts.
Let $\car\subset\cs$ be a full triangulated 
subcategory, closed under the formation of the coproducts in $\cs$ of
any set of its objects.
Form the category $\ct=\cs/\car$. It is easy to show that the
category $\ct$ contains all small coproducts, and that the 
natural map $\cs\ri\ct$ respects coproducts.
\entn

We have $\car\subset\cs$, with $\ct=\cs/\car$. The reader might 
imagine trying to apply Waldhausen's localisation theorem directly
to the triple $\car$, $\cs$ and $\ct$. There are two problems with this:
\be
\item Since $\car$, $\cs$ and $\ct$ contain all small coproducts of their
objects they tend to be huge categories. The classes of objects
are not small sets. This means that any permissible Waldhausen 
model would not be small, and there are set theoretic difficulties 
even in defining the Waldhausen \kth\ $K(\car)$, $K(\cs)$ and $K(\ct)$.

\item Even if we are willing to enlarge the universe and
define $K(\car)$ in this enlarged universe, we still get
nonsense. The fact that $\car$
contains countable coproducts of its objects permits us to
do the Eilenberg swindle, and show that $K(\car)$ is contractible.
Similarly for $K(\cs)$ and $K(\ct)$.
\ee

The useful way to produce a non-trivial, interesting example is
by passing to compact objects. The categories $\car$, $\cs$ and $\ct$
 each has a subcategory of compact objects.
The main theorem tells us:

\thm{My gen}
Let the notation be as in Notation~\ref{Notmainthm}.
 Assume further that there exist:
\be
\item A set of objects $S\subset \cs^c$, so that any subcategory
of $\cs$ containing $S$ and closed under triangles and coproducts is
all of $\cs$.
\item A set of objects $R\subset\car\cap \cs^c$,
so that any subcategory
of $\car$ containing $R$ and closed under triangles and coproducts is
all of $\car$.
\ee
Then the natural map $\car\ri\cs$ takes compact objects to compact
objects,
and so does the natural map $\cs\ri\ct$. 
In other words, we have a commutative 
diagram
\[
\CD
\car^c @>>> \cs^c
@>>> \ct^c\\
@VVV @VVV @VVV \\
\car @>>> \cs
@>>> \ct.
\endCD
\]
Of course the composite $\car^c\ri\cs^c\ri\ct^c$ must vanish, since
it is just the restriction to $\car^c$ of a vanishing functor on
$\car$.
We therefore have a factorisation of $\cs^c\ri\ct^c$ as
 \[
\CD
\cs^c @>>> \cs^c/\car^c
@>i>> \ct^c.
\endCD
\]
The functor $i\co\cs^c/\car^c
\ri\ct^c$ is fully faithful, and every object of $\ct^c$ is a direct
summand of an isomorph of an object in the image of $i$.
\ethm

\rmk{three proofs}
Thomason proved this theorem in the special case
where $\cs$ and $\ct$  are
the derived categories of quasi-coherent sheaves
on a scheme $X$ (respectively, on an open subset $U\subset X$).
In the generality above, the theorem may be found in
the first author's~\cite[Theorem~2.1]{Neeman92A}.
In a recent
 book~\cite{Neeman99} the first author
 generalises the theorem even further, to deal
 with the large cardinal case. There are now two proofs of
 Theorem~\ref{My gen}. The 
proof presented in the old paper~\cite{Neeman92A}, and the 
more general proof in
the book~\cite{Neeman99}. These two proofs
 are quite different from each other.
\ermk

\rmk{fully faith}
In the situation of Theorem~\ref{My gen} the map $\car\ri\cs$ is
fully faithful, and hence so is its restriction to $\car^c\ri\cs^c$.
Furthermore, every idempotent in the category $\car$ splits, because
$\car$ is closed under coproducts; 
see~\cite[Proposition~3.2 and Remark 3.3]{Bokstedt-Neeman93}.
Since every direct summand in $\car$ of a compact object is obviously
compact, every idempotent in $\car^c$ splits. It follows that
 $\car^c\subset\cs^c$ is closed under direct summands. 
\ermk

\section{The $\car$, $\cs$ and $\ct$ to which we
 apply Theorem~\protect{\ref{My gen}}}
\label{Smain}

In sections~\ref{Swaldh} and \ref{Sgenthom} we reviewed 
the general {\it K}--theoretic and
triangulated category results we will be using.
Now it is time to explain how we will apply them. We want to 
use the general theorems to deduce a \kth\ localisation theorem for
the Cohn localisation. 

The Cohn localisation begins with a ring $A$ and a set
$\s$ of morphisms $s_i\co P_i\ri Q_i$, as in Definitions~\ref{s-inve} 
and \ref{Dcohn}. To apply the results of 
Section~\ref{Sgenthom}, we need to choose suitable triangulated
categories $\car\subset\cs$, and $\ct=\cs/\car$. Our choices are:

\dfn{Fixed notation}
Let $A$ be a ring, $\s$ a set of maps of finitely generated,
projective $A$--modules. We define the triangulated categories
\be
\item $\cs=D(A)$ is the unbounded derived category of complexes
of $A$--modules.
\item We are given a set of maps $\s=\{s_i\co P_i\ri Q_i\}$. We can view these
as objects in $\cs=D(A)$ just by turning them into complexes
\[
\CD
\cdots @>>> 0 @>>> P_i @>s_i>> Q_i @>>> 0 @>>> \cdots
\endCD
\]
The category $\car\subset\cs$ is defined to be the smallest
triangulated subcategory of $\cs=D(A)$, which contains $\s$ and is
closed in $\cs$ under the formation of arbitrary coproducts of its
objects. 
\item $\ct$ is defined to be $\cs/\car$.
\ee
\edfn

\rmk{depen}
The categories $\car$, $\cs$ and $\ct$ depend on $A$ and on $\s$. In
most of this article we can view $A$ and $\s$ as fixed. For this reason
our notation makes no explicit mention of this dependence.
\ermk

Next we prove that 
our choices of $\car$, $\cs$ and $\ct$ satisfy the technical 
hypotheses of Theorem ~\ref{My gen}. We need a little lemma:

\lem{LSmain.0.5}
Let $A$ be a ring, and let $\cs=D(A)$. The object $A\in\cs$ is
the complex which is $A$ in degree $0$, and vanishes in all other
degrees. If $\cb\subset\cs$ is a triangulated subcategory 
which contains $A$ and is closed under coproducts,
then $\cb=\cs$.
\elem

{\bf Sketch of proof}\qua
This lemma is well-known and there are several 
proofs. We include a
sketch of one just for completeness.

Since $A\in\cb$ and $\cb$ is closed under direct sums, $\cb$ must contain
all free $A$ modules. Since $\cb$ is triangulated, it must contain
all bounded complexes of free $A$--modules. If $X\in\cb$ is a
bounded-above complex, then $X$ is quasi-isomorphic to a bounded
above complex $F$ of free modules. But $F$ is a direct limit of
its stupid truncations, all of which are bounded complexes of
free modules. The stupid truncations lie in $\cb$, and 
by~\cite[Remark 2.2]{Bokstedt-Neeman93} so does the direct limit 
$X\cong F$.

Now let $Y$ be an arbitrary (unbounded) object in $\cb$. Then $Y$ is
the direct limit of its (bounded above) \tstr\ truncations
$Y^{\leq i}$, all of which lie in $\cb$ by the above. Using
\cite[Remark 2.2]{Bokstedt-Neeman93} again, we conclude that $Y\in\cb$.
\endproof

In passing we mention the following corollary of 
Lemma~\ref{LSmain.0.5}:

\cor{Cmain.0.6}
As in Lemma~\ref{LSmain.0.5}, let $A$ be a ring and $\cs=D(A)$. 
An object $c\in\cs$ is compact if and only if it is isomorphic
to a perfect complex;
that is, if and only if $c$ is
isomorphic in
$\cs$ to a bounded chain complex of finitely generated, projective modules.
\ecor

\prf
In Example~\ref{Ecomp} we saw that every perfect complex 
is compact in $\cs$. We need to show that, up to isomorphism in
$D(A)$, these are the only
compact objects.

It is well known that the natural functor $D\perf(A)\ri D(A)$ is
fully faithful; $D\perf(A)$ may be viewed as a full subcategory of 
$\cs=D(A)$. Let $\widetilde\car$ be the subcategory of $D(A)$ 
containing all objects isomorphic to objects of 
$D\perf(A)\subset D(A)$.
The subcategory $\widetilde\car$ is triangulated, and from
\cite[Proposition 3.4]{Bokstedt-Neeman93} any 
direct summand of an object in $\widetilde\car$ lies in
$\widetilde\car$. The subcategory
$\widetilde\car$ also contains the object $A\in \cs$, which generates $\cs$ 
by Lemma~\ref{LSmain.0.5}. By
\cite[Lemma 2.2]{Neeman92A} it follows that $\widetilde\car$ contains all
the compact objects.
\eprf

\pro{LSmain.1}
Let the triangulated categories $\car\subset\cs$, $\ct=\cs/\car$
be as in Definition~\ref{Fixed notation}. Then
$\car\subset\cs$, $\ct=\cs/\car$ satisfy the
hypotheses of Theorem~\ref{My gen}.
\epro

\prf
The category $\cs=D(A)$ clearly contains coproducts of its objects,
and by its definition $\car\subset\cs$ is closed in $\cs$ under 
coproducts. That is, the notation is as in Notation~\ref{Notmainthm}.

It remains to verify the hypotheses~\ref{My gen}(i) and (ii). 
For $S\subset\cs^c$ take the set $\{A\}$. Lemma~\ref{LSmain.0.5}
tells us that \ref{My gen}(i) holds. For the set $R\subset\car$
of \ref{My gen}(ii) we take $\s$. The definition of $\car$
is as the smallest triangulated subcategory of $\cs$, closed under
coproducts and containing $\s$. To prove \ref{My gen}(ii), it
suffices to establish that $\s\subset\cs^c$.

But every object of $\s$ is a chain complex
\[
\CD
\cdots @>>> 0 @>>> P_i @>s_i>> Q_i @>>> 0 @>>> \cdots
\endCD
\]
with $P_i$ and $Q_i$ finitely generated and projective. By
Corollary~\ref{Cmain.0.6} (or even by Example~\ref{Ecomp})
it follows that every object in $\s$ is compact in $\cs=D(A)$.
\eprf

Since the hypotheses of Theorem~\ref{My gen} hold, so does its
conclusion. We deduce a diagram of triangulated categories:
\[
\xymatrix{
\car^c \ar[rr] & & \cs^c \ar[dr]
\ar[rr]^-{\pi} & &
\ct^c    \\
 & & & \cs^c/\car^c \ar[ur]_{i} & 
}
\] 
\rmk{Rchoices}
Now we make our choices of permissible Waldhausen categories.
We let $\bS=C\perf(A)$ be the Waldhausen category of all
perfect chain complexes of $A$--modules, with morphisms the
chain maps, weak equivalences the homotopy equivalences, and
cofibrations the degreewise split monomorphisms. 
Clearly $\bS=C\perf(A)$ is a model for $D\perf(A)$.
As in 
Lemma~\ref{cheapmodels} we produce $\bR$ and $\bSR$.
The 
category $\bR$ is the full subcategory 
of all objects in $\bS=C\perf(A)$
which become isomorphic in $D(A)$ to objects in $\car^c$.
The category $\bSR$ has the same objects,
morphisms and cofibrations as $\bS$, but the weak equivalences
are any morphisms whose mapping cones lie in $\bR$. 

Slightly more delicate is our choice for $\bT$. For this
we need:

\dfn{Dcatmod}
For any ring $B$, the category $C(B,\aleph_0)$ will be defined
as follows. The objects are certain chain complexes of projective
$A$--modules, to be specified below. The morphisms are the chain maps,
the weak equivalences are the quasi-isomorphisms, and the cofibrations
are the degreewise split monomorphisms. The restrictions on the
objects are given by specifying that $C(B,\aleph_0)$ is the
smallest category which:
\be
\item Contains the perfect complexes.
\item Is closed under countable direct sums.
\item Is closed under the formation of mapping cones.
\ee
If $A$ is not just any ring, but comes with a set $\s$ of maps of finitely
generated projective $A$--modules, then 
$C(A,\s,\aleph_0)$ has the same objects, morphisms and cofibrations
as $C(A,\aleph_0)$. Only the weak equivalences change. A morphism
in $C(A,\s,\aleph_0)$ is a weak equivalence if its mapping cone
maps to $\car\subset\cs$ under the functor $C(A,\aleph_0)\ri
D(A)=\cs$.
\edfn

\nin
Our choice for $\bT$ is to take the full Waldhausen subcategory 
of $C(A,\s,\aleph_0)$ whose objects are isomorphic in
$\ct=D(A)/\car$ to compact objects.\footnote{The
reason for this definition is that $\bT$ must be
essentially small to have a \kth. Hence we only include countable
coproducts.}
\ermk

\lem{LOKperm}
If we take the diagram of permissible Waldhausen categories
\[
\xymatrix{
\bR \ar[rr] & & \bS \ar[dr]
\ar[rr]^-{\pi} & &
\bT    \\
 & & & \bSR \ar[ur]_{i} & 
}
\] 
and pass to derived categories, we obtain (up to equivalence):
\[
\xymatrix{
\car^c \ar[rr] & & \cs^c \ar[dr]
\ar[rr]^-{\pi} & &
\ct^c    \\
 & & & \cs^c/\car^c \ar[ur]_{i} & 
}
\] 
\elem

\prf
By Corollary~\ref{Cmain.0.6}, there is an equivalence of categories
$D\perf(A)\ri \cs^c$. This makes $\bS=C\perf(A)$
a model for $\cs^c$. The definitions of the categories $\bR$
and $\bSR$, together with
Lemma~\ref{cheapmodels},
 make the lemma immediate for the  
part of the diagram
\[
\xymatrix{
\bR \ar[rr] & & \bS \ar[dr]
 & &
    \\
 & & & \bSR  & 
}
\] 
The slight subtlety comes from $\bT$. 
The key point is that the derived
category $D(A,\s,\aleph_0)$ of $C(A,\s,\aleph_0)$
maps fully faithfully to $\ct=\cs/\car$.
This is not
entirely trivial;\footnote{This is the only point in the article where 
we appeal to a theorem about triangulated categories which cannot be
found in the three basic 
references~\cite{BeiBerDel82,Bokstedt-Neeman93,Neeman92A}.}
it may be found in
\cite[Proposition~4.4.1]{Neeman99}. Since  $C(A,\s,\aleph_0)$ is
closed under countable direct sums so is 
$D(A,\s,\aleph_0)$.
 From~\cite[Remark 3.3]{Bokstedt-Neeman93}
we conclude that any idempotent in 
$D(A,\s,\aleph_0)$ splits. 
But $D(A,\s,\aleph_0)$ is equivalent to a full subcategory of
$\ct$ containing the image of $D\perf(A)$ and closed under
direct summands. It follows that
$\ct^c\cap D(A,\s,\aleph_0)$ is equivalent to $\ct^c$.
Since $\bT$ is defined to be the preimage in $C(A,\s,\aleph_0)$
of
$\ct^c\cap D(A,\s,\aleph_0)$, its derived 
category must be equivalent to $\ct^c$ (again by Lemma~\ref{cheapmodels}).
\eprf

\cor{Cmain.2}
In the sequence
\[
\CD
K(\bR) @>>> K(\bS) @>{K(\pi)}>> K(\bT)
\endCD
\]
the spectrum $K(\bR)$ is the $(-1)$--connected cover of the
homotopy fiber of $K(\pi)$. 
Furthermore, $K(\bS)$ agrees with $K(A)$, the
Quillen \kth\ of $A$.
\ecor

\prf
The statement about the fiber is immediate from 
Theorems~\ref{Local} and \ref{Graycof}.
By definition $\bS=C\perf(A)$, and hence $K(\bS)=K\big(C\perf(A)\big)$.
The assertion $K\big(C\perf(A)\big)=K(A)$
may be found in Gillet's~\cite[6.2]{Gillet81}.
\eprf

\section{The map $\ct\ri D(\rs)$}
\label{TtoB}

Let $A$ be an associative ring, and let $\s$ be a set of maps of finitely
generated, projective $A$--modules. Let the categories $\cs=D(A)$, 
$\car\subset\cs$ and $\ct=\cs/\car$ 
be as in Definition~\ref{Fixed notation}.

In the previous section we showed that in the sequence
\[
\CD
K(\bR) @>>> K(A) @>{K(\pi)}>> K(\bT)
\endCD
\]
$K(\bR)$ is the $(-1)$--connected cover of the homotopy
fiber of the map ${K(\pi)}$.
Next we want necessary and sufficient conditions 
for $\ct^c$ to be $D\perf(\rs)$. Under these conditions
\[
K(\bT)=K\big(C\perf(\rs)\big)=K(\rs)
\]
where the last equality is by Gillet's~\cite[6.2]{Gillet81}.
The sequence above becomes, up to some nonsense in 
the $(-1)$--homotopy groups, a homotopy fibration
\[
\CD
K(\bR) @>>> K(A) @>>> K(\rs)
\endCD
\]
and this is what we are after. 

The first step is to find a functor comparing $\ct^c$ and $D\perf(\rs)$.
We define it at the level of unbounded complexes.
Let us remind the reader first of the tensor product of unbounded
complexes.

\rmd{RT} Let $B$ be any $(A-A)$--bimodule.
The {\it derived tensor product with $B$} is 
a triangulated, coproduct-preserving functor
$D(A)\ri D(A)$. We will denote it
\[
X\mapsto B\oti X.
\]
If $B$ is an $A$--algebra, we can view this as a functor 
$D(A)\ri D(B)$.
The existence of this functor was first proved by 
Spaltenstein~\cite{Spaltenstein88}. A very short proof of
the existence may be found in~\cite[Theorem~2.14]{Bokstedt-Neeman93}.
Very concretely, to form $B\oti X$ we take 
a {\it K}--projective resolution $P\ri X$, and define
\[
B\oti X\eq B\otimes_A^{}P.
\]
Recall
that a map $P\ri X$ is a $K$--projective resolution
if
\be
\item $P\ri X$ is a quasi-isomorphism.
\item Any chain map $P\ri Y$, with $Y$ acyclic, is null homotopic.
\ee
The references above prove the existence of $K$--projective resolutions%
\footnote{The terminology is a little misleading. Note that a
$K$--projective resolution is not necessarily a projective resolution
in the ordinary sense of the word.  Any null homotopic complex $N$
is a $K$--projective resolution of the zero complex. But $N$ is not
necessarily a complex of projectives.}.

In this article, we consider tensor products both in the
category of modules and in the derived category.  We try to be careful
to distinguish them in the notation.  
\ermd

It will be helpful to note that, for the
categories $\bT\subset 
C(A,\s,\aleph_0)$ of Remark~\ref{Rchoices} and
Definition~\ref{Dcatmod}, the ordinary 
tensor product agrees with the derived
tensor product.

\lem{K-projective}
The objects $P\in C(A,\aleph_0)$ are all $K$--projective. 
\elem

\prf
The perfect complexes are clearly $K$--projective. Furthermore any 
coproduct of $K$--projectives is $K$--projective, and any mapping
cone on a map of $K$--projectives is $K$--projective.
\eprf

\lem{LTtoD.1}
Let $A\ri B$ be any $\s$--inverting ring homomorphism 
(Definition~\ref{s-inve}). By Reminder~\ref{RT} there is a functor
$D(A)\ri D(B)$, taking $X\in D(A)$ to $B\oti X\in D(B)$. We
assert that this functor factors uniquely as
\[
\CD
D(A)=\cs @>\pi>> \ct @>T>> D(B)
\endCD
\]
where $\pi\co\cs\ri\ct=\cs/\car$ is as in 
Definition~\ref{Fixed notation}, and $T$ respects 
coproducts. Furthermore, the
functor $T\co\ct\ri D(B)$ takes compact objects to compact
objects. 
\elem

\prf
Let $s_i\co P_i\ri Q_i$ be a map in $\s$. Tensoring
with $B$ takes it to an isomorphism. Hence tensoring
with $B$ takes the chain complex
\[
\CD
\cdots@>>>  0 @>>>  P_i @>s_i>>  Q_i @>>> 0 @>>>\cdots
\endCD
\]
to an acyclic complex.  Therefore the functor 
$X\mapsto B\oti X\co D(A) \ri D(B)$ kills all
the objects in $\s$.  Since derived tensor product preserves triangles
and coproducts, the subcategory of $\cs$ annihilated by 
$X\mapsto B\oti X$ must be
closed under triangles and coproducts, and therefore contains all of
$\car$.  By the universal property of the Verdier quotient
$\ct=\cs/\car$, there is a unique factorisation
\[
\CD
D(A)=\cs @>\pi>> \ct @>T>> D(B)
\endCD
\]
and since $T\pi$ and $\pi$ respect coproducts, so does $T$.
It remains to show that $T$ takes $\ct^c\subset\ct$ to 
${\{D(B)\}}^c\subset D(B)$.

It is clear that the map $T\pi\co\cs\ri D(B)$ takes a bounded
complex of finitely generated projective $A$--modules 
to a bounded complex of 
finitely generated projective
$B$--modules; the map just tensors with $B$.  
By Corollary~\ref{Cmain.0.6}, this says that the 
functor $T\pi$ takes $\cs^c$ to ${\{D(B)\}}^c\subset D(B)$.
The last statement of Theorem~\ref{My gen} says that every object
$t\in\ct^c$ is a direct summand of $\pi(s)$, with $s\in\cs^c$.
Hence $T(t)$ is a direct summand of the compact $T\pi(s)$,
and must therefore be compact.
\eprf

\smr{Summary2}
For the special
$\car$, $\cs$ and $\ct$ of Definition~\ref{Fixed notation},
Lemma~\ref{LTtoD.1} allows us
 to extend the diagram of Theorem~\ref{My gen} to:
$$\xymatrix{
\car^c \ar[rr] \ar[dd]& & D\perf(A) \ar[dr]
  \ar[rr]^-{\displaystyle{\pi}}\ar[dd] & &
  \ct^c  \ar[dd]\ar[rr]^{T} & & D\perf(B)\ar[dd]& \\
& & & {D\perf(A)}/{\car^c} \ar[ur]_{i} & & & &
  \!\!\!\!\!\!\!\!\!\!\!\!(**)\\
\car\ar[rr] &
  & D(A) \ar[rr]_-{\displaystyle{\pi}} &
   & \ct\ar[rr]_{T} & & D(B) &}$$
Now we construct $\bU$, a permissible Waldhausen model for $D\perf(B)$. 
First consider the category $C(B,\aleph_0)$ of Definition~\ref{Dcatmod}.
The category $\bU$ is defined to be the full
subcategory of $C(B,\aleph_0)$ whose objects are quasi-isomorphic to
perfect complexes. The map $X\mapsto B\otimes_A X$ defines an exact
functor $C(A,\s,\aleph_0)\ri C(B,\aleph_0)$, which takes $\bT\subset
C(A,\s,\aleph_0)$ to $\bU\subset C(B,\aleph_0)$.
The diagram of permissible Waldhausen categories
$$\xymatrix{
\qquad&\bR\ar[rr] & & \bS \ar[dr]
\ar[rr]^-{\displaystyle{\pi}} & &
 \bT  \ar[rr]^{T} & & \bU&\qquad \\
& & & & \bSR \ar[ur]_{i} & & & &
}$$
gives, when we pass to derived categories, precisely the top
row of $(**)$. 

Before applying Waldhausen's \kth\ to
this diagram, it is helpful to extend it a little bit. 
Put $\bD=C\perf(B)$.
Remember that
tensor product with $B$ takes perfect complexes to perfect complexes.
Therefore we have a commutative square, where the vertical maps
are induced by tensor product with $B$:
\[
\CD
C\perf(A)@= \bS @>>> \bT @.\sub @. C(A,\s,\aleph_0)\\
@VVV      @VVV   @VVV @.    @VVV \\
C\perf(B)@= \bD @>>> \bU @.\sub @. C(B,\aleph_0)
\endCD
\]
The map from $C\perf(A)= \bS$ to $C\perf(B)= \bD$ clearly factors
through $\bSR$; after all, the only change from $\bS$ to $\bSR$
is in the weak equivalences, and the larger class of weak equivalences
in $\bSR$ maps to weak equivalences in $\bD$. We therefore
have a commutative diagram:
\[
\xymatrix{
\qquad&\bR\ar[rr] & & \bS \ar[dr]
  \ar[rr]^-{\displaystyle{\pi}} & &
  \bT  \ar[rr]^{T} & & \bU&\qquad (*\!*\!*) \\
& & & & \bSR \ar[ur]_{i}\ar[rr] & &\bD\ar[ur] & &
}
\]
In this diagram, the map
$\bD\ri\bU$ induces an equivalence of derived categories; both
$\bD$ and $\bU$ are models for $D\perf(B)$. By Theorem~\ref{Approx}
the map $K(\bD)\ri K(\bU)$ is a homotopy equivalence.
Gillet's~\cite[6.2]{Gillet81} tells us that 
\[
K(\bS)=K(C\perf(A))=K(A), \qquad\hbox{ and }\qquad
K(\bD)=K(C\perf(B))=K(B).
\]
In \kth, $(*\!*\!*)$ yields:
 \[
\xymatrix{
K(\bR) \ar[rr] & & K(A) \ar[dr]
\ar[rr]^-{\displaystyle{K(\pi)}} & &
 K(\bT)  \ar[rr]^{K(T)} & & K(\bU) \\
 & & & K(\bSR) \ar[ur]_{K(i)}\ar[rr] & & K(B)\ar[ur]_{\simeq} & 
}
\]
In the remainder of this section we will analyse necessary and
sufficient conditions for the functor $T\co\ct^c\ri D\perf(B)$ 
to be an equivalence
of categories. When it is, it follows from  Theorem~\ref{Approx}
that $K(T)\co K(\bT)\ri K(\bU)$ is a homotopy equivalence. From
the diagram and our previous discussion it then follows that, up
to nonsense in degree $(-1)$, $K(\bR)$ is the homotopy fiber 
of the natural map $K(A)\ri K(B)$.
\esmr

\lem{LTtoD.2}
Let the notation be as in Lemma~\ref{LTtoD.1}. If the functor
$T\co\ct^c\ri D\perf(B)$ is an equivalence, then
\be
\item The functor $\pi\co D(A)=\cs\ri\ct$
induces a homomorphism
\[
\CD
\Hom_{D(A)}^{}(A, A) @>>> \ct(\pi A,\pi A)
\endCD
\]
which can be naturally identified with
$A\op\ri B\op$.
\item For all $n\neq 0$, $\ct(\pi A,\T^n\pi A)=0$.
\ee
\elem

\prf
If $T$ is an equivalence, then for all $s,t\in\ct^c\subset\ct$ we must have
\[
\ct(s,t)=\Hom_{D(B)}^{}(Ts,Tt).
\]
Put $s=\pi A$ and $t=\T^n\pi A$. Then $Ts=T\pi A=B\otimes_A^{}A=B$,
and $Tt=\T^n B$. This gives
\[
\ct(\pi A,\T^n\pi A)=\Hom_{D(B)}^{}(B,\T^n B)
\]
and the right hand side is $B\op$ in degree zero, and vanishes if $n\neq0$.
Moreover, the induced map $A\op\ri B\op$ is the natural homomorphism.
\eprf
 
The most interesting case is $B=\rs$. The proof of the next
Proposition is a small modification of ideas that may be
found in
Rickard's~\cite{Rickard89b}. 

\pro{PTtoD.3}
Let the notation be as in Lemma~\ref{LTtoD.1}, but 
with $B=\rs$. The following are
equivalent:
\be
\item The functor $T\co\ct\ri D(\rs)$ is an equivalence.
\item The restriction to compact objects, that is
 $T\co\ct^c\ri D\perf(\rs)$, is an equivalence.
\item
$\ct(\pi A,\T^n\pi A)=\left\{
\begin{array}{lll}
\rrs &\quad & \hbox{if }n=0 \\
\,0 &\quad & \hbox{if }n\neq0
\end{array}
\right.$\par
where the isomorphism $\ct(\pi A,\pi A)=\rrs$ is as
$A\op$--algebras.
\ee
\epro

\prf
(i)$\Longrightarrow$(ii) is clear; if $T$ is an equivalence, 
then it restricts to an equivalence on compact objects. 
(ii)$\Longrightarrow$(iii) was proved in Lemma~\ref{LTtoD.2}. It
remains to prove (iii)$\Longrightarrow$(i).

Assume now that (iii) holds. Among other things, we know that
there
is some isomorphism $\ct(\pi A,\T^n\pi A)=\rrs$ of $A\op$--algebras.
We first want to 
show how it follows that the functor $T$ induces an isomorphism
\[
\CD
\ct(\pi A,\pi A)
@>T>>\Hom_{D(\rs)}^{}(\rs,\rs)=\rs.
\endCD
\]
We have ring homomorphisms
\[
\CD
\Hom_{D(A)}^{}(A,A) @>\pi>>
\ct(\pi A,\pi A)
@>T>>\Hom_{D(\rs)}^{}(\rs,\rs)
\endCD.
\]
That is,
\[
\CD
A\op @>\pi>> \ct(\pi A,\pi A) @>T>> \rrs.
\endCD
\]
The composite $T\pi\co A\op\ri\rrs$ is easily computed to be the natural
map. By hypothesis (iii), there is an isomorphism 
$\ct(\pi A,\pi A)\cong\rrs$, as $A\op$--algebras. But then both
$\pi\co A\ri {\ct(\pi A,\pi A)}\op$ and $T\pi:A\ri \rs$ 
are initial in the category of $\s$--inverting
ring homomorphisms. Hence the map $T\co \ct(\pi A,\pi A)\ri\rrs$
must be an isomorphism.

Let $\cc\subset\ct$ be the full subcategory with objects
\[
\text{\rm Ob}(\cc)=
\left\{x\in {\rm Ob}(\ct)\,\left|\begin{array}{c} 
\forall n\in\zz\text{ the map induced by }T\\
\ct(\pi A,\T^n x)\ri\Hom_{D(\rs)}^{}(T \pi A,\T^n Tx)\\
\text{is an isomorphism}\end{array}\right.\right\}.
\]
Since $\pi A\in\ct$ and $T\pi A=\rs\in D(\rs)$ are both
compact, the category $\cc$ is closed under
coproducts. It is clearly closed under triangles, and by (iii) and
the above it
contains $\pi A$. Its inverse image under the projection map
\[
\CD
\pi\co D(A)=\cs @>>> \cs/\car=\ct
\endCD
\]
is a triangulated subcategory, closed under coproducts and containing 
$A$. Now Lemma~\ref{LSmain.0.5} tells us that $\pi^{-1}\cc=\cs$, and
hence $\cc=\ct$.

Next let $\cd\subset\ct$ be the full subcategory with objects
\[
\text{\rm Ob}(\cd)=
\left\{x\in {\rm Ob}(\ct)\,\left|\begin{array}{c} 
\forall n\in\zz\text{ and }\forall y\in\ct \text{ the map}\\
\ct(x,\T^n y)\ri\Hom_{D(\rs)}^{}(T x,\T^n Ty)\\
\text{is an isomorphism}\end{array}\right.\right\}.
\]
By the above, $\cd$ contains $\pi A$. It is clear that $\cd$ is closed
under triangles and coproducts. As above, it follows that $\pi^{-1}\cd=\cs$,
and hence $\cd=\ct$. 

This proves that $T$ is fully faithful. It embeds $\ct$ as a full,
triangulated subcategory of $D(\rs)$, closed under coproducts and containing
$T\pi A=\rs$. Applying Lemma~\ref{LSmain.0.5} to the inclusion
$\ct\subset D(\rs)$, we conclude it must be an equivalence.
\eprf

\section{The case $n>0$}
\label{Sproofs}

Let the notation be as in Proposition~\ref{PTtoD.3}. That is, $A$ is a ring,
$\s$ is a set of maps of finitely generated, projective 
$A$--modules, $\rs$ is
the Cohn localisation, $\car$, $\cs$ and $\ct$ are the triangulated categories
of Definition~\ref{Fixed notation}, and $T\co \ct\ri D(\rs)$ is the
functor of Lemma~\ref{LTtoD.1}. Proposition~\ref{PTtoD.3} tells us
that everything is reduced to computing the groups $\ct(\pi A,\T^n\pi A)$.
In the next three sections we will prove
\be
\item If $n>0$, then $\ct(\pi A,\T^n\pi A)=0$. [This section].
\item If $n=0$, then the ring homomorphism 
$A\op\ri \ct(\pi A,\T^n\pi A)$ agrees
with $A\op\ri\rrs$. [Section~\ref{n=0}].
\item The groups $\ct(\pi A,\T^n\pi A)$ vanish for all $n<0$ if and
only if the groups $\Tor_n^A(\rs,\rs)$ vanish for all $n>0$.
[Section~\ref{n<0}].
\ee

In every paper there comes a time for hard work. The day of reckoning has
come in this paper. We now have to prove something. The key tool in the 
proofs is:

\pro{j^*}
The functor $\pi\co\cs\ri\ct$ has a right adjoint $G\co\ct\ri\cs$.
The unit of adjunction $\eta_x^{}\co x\ri G\pi x$ can be completed
to a distinguished triangle
\[
\CD
k @>>> x @>\eta_x^{}>> G\pi x @>>> \T k.
\endCD
\]
In this triangle, the object $k$ lies in $\car$.
\epro

\prf
See \cite[Lemma~1.7]{Neeman92A}. The notation there is slightly different; 
the functor we have been calling $\pi\co\cs\ri\ct$ is called $j^*$ there,
and the adjoint we call $G$ goes by the name $j_*$ there.
\eprf

\rmd{Rj^*}
An object $y\in \cs$ is called {\em $\s$--local} (or just
{\em local} if $\s$ is understood) if for all
$r\in\car$ we have $\cs(r,y)=0$. 
If $x$ is any object of $\cs$, then
the object $G\pi x$ is 
local; after all 
\[
\begin{array} {rclcl}
\cs(r,G\pi x)&=& \ct(\pi r, \pi x)&\qquad&\hbox{by adjunction}\\
&=& \ct(0, \pi x)&\qquad&\hbox{since }r\in\car,\hbox{ hence }\pi r=0 \\
&=& 0. & & 
\end{array}
\]
\ermd

\lem{Gtrunc}
If $x\in\cs$ is a local object, then so are its \tstr\ truncations
$x^{\leq n}_{}$ and $x^{\geq n}_{}$.
\elem

\prf
Pick a $\s$--local object $x$ and an integer $n\in\zz$;
we will show first that $x^{\geq n}$ is $\s$--local.
Without loss of generality we may assume $n=0$.
To prove that $x^{\geq 0}$ is $\s$--local,
take any $s\in\s$. We will show that $\cs(s,x^{\geq 0})=0$.
Assume for a second that we have shown this, for all $s\in\s$.
This will mean that the full subcategory $\cc\subset\cs$ 
whose objects are
\[
\text{\rm Ob}(\cc)=\{
c\in\cs\mid \forall n\in\zz,\,\,\,\cs(\T^n c,x^{\geq 0})=0
\}
\]
contains $\s$. But $\cc$ is clearly triangulated and closed under
coproducts. Hence $\car\subset\cc$, which means that $x^{\geq 0}$
is local.

Hence it needs to be shown that, for any $s\in\s$ and 
any local $x$, $\cs(s,x^{\geq 0})=0$. Let $s$ be the chain complex
\[
\CD
\cdots @>>> 0 @>>>  P_i @>s_i>> Q_i @>>> 0 @>>>\cdots 
\endCD
\]
Now $P_i$ is in some degree $m$ and $Q_i$ is in degree $m+1$.
There are two cases:

\nin
{\bf Case 1}\qua
If $m\leq -2$, then $s$ is a complex concentrated 
in degrees $\leq -1$; that is, $s\in\cs^{\leq -1}$. 
But $x^{\geq 0}$ is in $\cs^{\geq0}$. 
Hence all maps $s\ri x^{\geq 0}$
vanish.

\nin
{\bf Case 2}\qua
Suppose $m\geq -1$.
The \tstr\ gives a distinguished triangle
\[
\CD
x^{\leq -1} @>>> x @>>> x^{\geq 0} @>w>> \Sigma x^{\leq -1}.
\endCD
\]
For any map $s\ri x^{\geq 0}$, the composite 
\[
\CD
s @>>>  x^{\geq 0}@>w>> \Sigma x^{\leq -1}
\endCD
\]
is a map from a bounded above complex of projectives
$s$ to some object in $\cs=D(A)$, and hence
it is represented by a chain map. But the chain
complex $s$ lives in degrees $m$ and
$m+1$, both of which are $\geq-1$, while the complex
$\Sigma x^{\leq -1}$ lies in $\cs^{\leq-2}$. Hence the map vanishes.
 From the triangle we deduce that the map
$s\ri x^{\geq 0}$ must factor as
\[
\CD
s @>>> x @>>> x^{\geq 0}.
\endCD
\]
Now $x$ is $\s$--local by hypothesis, and hence 
any map $s \ri x$ vanishes.

This proves that $x^{\geq0}$ is $\s$--local. We have a triangle
\[
\CD
x^{\leq -1} @>>> x @>>> x^{\geq 0} @>w>> \Sigma x^{\leq -1}~.
\endCD
\]
The long exact sequence for $\cs(r,-)$, with $r\in\car$,
allows us 
to deduce that
$x^{\leq -1}$ is also $\s$--local. Shifting by powers of $\T$,
we deduce that $x^{\leq n}$ is $\s$--local for any $n\in\zz$.
\eprf

\lem{L6.1}
Let the notation be as above.
If $x\in\cs^{\leq n}$, 
then so is $G\pi x$.
\elem

\prf
We may assume without loss that $n=0$. Pick any 
$x\in\cs^{\leq 0}$. By Reminder~\ref{Rj^*},
$G\pi x$ is local.  By Lemma~\ref{Gtrunc}
so is ${\{G\pi x\}}^{\leq 0}$.

Now the unit of adjunction $\eta_x\co x\ri G\pi x$
is a map from an object $x\in\cs^{\leq 0}$, and
therefore factors (uniquely) as
\[
\CD
x @>\alpha>> {\{G\pi x\}}^{\leq0} @>f>> G\pi x~.
\endCD
\]
On the other hand, we have a triangle
\[
\CD
k @>>> x @>\eta_x>> G\pi x @>>> \T k
\endCD
\]
with $k\in\car$. The composite
\[
\CD
k @>>> x @>\alpha >> {\{G\pi x\}}^{\leq0}
\endCD
\]
must vanish, since $k\in \car$ and ${\{G\pi x\}}^{\leq 0}$
is local. It follows that $\alpha$ factors as
\[
\CD
x @>\eta_x>> G\pi x @>g>> {\{G\pi x\}}^{\leq0}.
\endCD
\]
The composite
\[
\CD
x @>\eta_x>> G\pi x @>g>> {\{G\pi x\}}^{\leq0} @>f>> G\pi x
\endCD
\]
is $\eta_x$, by construction of $f$ and $g$.
It follows that
\[
\CD
x @>\eta_x>> G\pi x @>1-fg>> G\pi x
\endCD
\]
vanishes. From the triangle
\[
\CD
k @>>> x @>\eta_x>> G\pi x @>>> \T k
\endCD
\]
we have that $1-fg$ must factor through a map $\T k\ri G\pi x$.
But as $\T k\in\car$ and $G\pi x$ is local, $1-fg$ must vanish.
In other words, $fg=1$.

But then the identity on $G\pi x$ factors through an object
in $\cs^{\leq 0}$. Therefore $1\co H^n(G\pi x)\ri H^n(G\pi x)$ 
vanishes for all $n>0$. This means that for $n>0$ we have 
$H^n(G\pi x)=0$. In other words, $G\pi x\in\cs^{\leq0}$.
\eprf

\cor{Cvan}
Let the notation be as above. 
For any $n>0$ we have 
\[
\ct(\pi A,\T^n\pi A)=0.
\]
\ecor
\proof
We know $A\in\cs^{\leq0}$, and from Lemma~\ref{L6.1} we deduce
$G\pi A\in\cs^{\leq0}$. Now we compute
\[
\begin{array}{rclclr}
\ct(\pi A,\T^n\pi A) &=& \cs(A,\T^n G\pi A) &\quad&\text{by adjunction} & \\
&=& H^n(G\pi A)&\quad&\text{because }\cs(A,\T^n X)=H^n(X) & \\
&=& 0&\quad&\text{because }G\pi A\in\cs^{\leq0}. & \qquad~\Box
\end{array}
\]

\section{The case $n=0$}
\label{n=0}

In this section we will show that $\ct(\pi A,\pi A)=\rrs$. 
We know that $\ct(\pi A,\pi A)$ is a ring, and comes with a natural
ring homomorphism
\[
\CD
A\op=\cs(A,A) @>>> \ct(\pi A,\pi A)=B\op.
\endCD
\]
What we prove is that the ring homomorphism $A\ri B$ above is
the initial $\s$--inverting homomorphism.

\lem{Ln=0.1}
The unit of adjunction gives us a map $\eta_A^{}\co A\ri G\pi A$.
Applying the functor $H^0$ gives a map
\[
\CD
A=H^0(A) @>H^0(\eta_A^{})>> H^0(G\pi A).
\endCD
\]
We assert that this agrees with the natural homomorphism
\[
\CD
A @>>> \ct(\pi A,\pi A)\op=B.
\endCD
\]
That is, there is an isomorphism of left $A$--modules
$B\ri H^0(G\pi A)$, commuting with the inclusion of $A$.
\elem

\proof
It is clear how to give an isomorphism of sets
$H^0(G\pi A)\cong
\ct(\pi A,\pi A)$. We have
\[
\begin{array}{rclcl}
\ct(\pi A,\pi A) &=& \cs(A,G\pi A) &\qquad & \text{by adjunction} \\
&=& H^0(G\pi A). &\qquad & 
\end{array}
\]
As sets, we have an equality $\ct(\pi A,\pi A)= \ct(\pi A,\pi A)\op$. 
Let us call the isomorphism of sets given above
$\p\co \ct(\pi A,\pi A)\op\ri H^0(G\pi A)$.
We have a triangle 
\[
\xymatrix@R-10pt@C+5pt{
  & & \ct(\pi A,\pi A)\op\ar^{\p}[dd] \\
A\,\, \ar[rru]\ar_{H^0(\eta_A^{})}[rrd] & & \\
 & & H^0(G\pi A).
}
\]
We need to show that the triangle commutes, and that $\p$ is
a map of left $A$--modules. For a second let us suppose we know that
$\p$ is a homomorphism of left $A$--modules. Then all three
maps are $A$--module homomorphisms, and the commutativity can be checked 
by evaluating the maps on the single element $1\in A$.
We leave this to the reader.

It remains to show that $\p$ is a homomorphism of left
$A$--modules. We must show that $\p$ 
takes right multiplication by $a\in A$ in
the ring $\ct(\pi A,\pi A)$ 
to left multiplication by $a$ in $H^0(G\pi A)$.

Therefore we let $x$ be any element of $\ct(\pi A,\pi A)$, and
let $a\in A$. Then $xa$ is an element of the ring
$\ct(\pi A,\pi A)$;
it is the composite
\[
\CD
\pi A @>\pi \rho_a>> \pi A @>x>> \pi A,
\endCD
\]
where $\rho_a\co A\ri A$ is the map induced by right multiplication
by $a$. The naturality of $\eta\co 1\ri G\pi$ gives a commutative
diagram
\[
\CD
A @>\rho_a>> A @. \\
@V{\eta_A^{}}VV  @VV{\eta_A^{}}V @. \\
G\pi A @>>{G\pi\rho_a }> G\pi A @>>{Gx}> G\pi A.
\endCD
\]
Consider the image of $1\in A$. We have
\[
\begin{array}{rclclr}
\p(xa) &=&G(x\,\,\pi\rho_a^{})\,\eta_A^{}(1) &\qquad & 
\text{adjunction formula} & \\
  &=& Gx\,\, G\pi\rho_a^{}\,\,\eta_A^{}(1) &\qquad & G 
\text{ respects composition} & \\
&=& Gx\,\, \eta_A^{}\,\, \rho_a(1)&\qquad & 
\text{by commutative diagram above} & \\
&=& Gx\,\, \eta_A^{} (a)&\qquad & \text{since }\rho_a(1)=a & \\
&=& a\, Gx\,\, \eta_A^{} (1)&\qquad & \text{since }Gx\,\eta_A^{}
\text{ is a homomorphism} & \\
& & & & \text{of left }A\text{--modules} & \\
&=&a\,\p(x). & & & \qquad~\Box
\end{array}
\]

\lem{Ln=0.2}
The ring homomorphism $A\ri \ct(\pi A,\pi A)\op$ is $\s$--inverting.
\elem

\prf
Let $s_i\co P_i\ri Q_i$ be any element
of $\s$. Because $G\pi A$ is local we know that,
for any $n\in\zz$, $\cs(\T^ns_i,G\pi A)=0$.
 From the distinguished triangle
\[
\CD
P_i@>>> Q_i @>>>  s_i @>>> \T P_i
\endCD
\]
we conclude that the natural map
\[
\CD
\cs(Q_i, G\pi A)@>>> \cs (P_i, G\pi A) 
\endCD
\]
is an isomorphism. But both $P_i$ and $Q_i$ are projective
$A$--modules, viewed as complexes concentrated in degree $0$.
Therefore the isomorphism above is the natural map
\[
\CD
\Hom_A^{}\big(Q_i, H^0(G\pi A)\big)@>>> 
\Hom_A^{}\big(P_i, H^0(G\pi A)\big). 
\endCD
\]
Put $B=\ct(\pi A,\pi A)\op$. As left $A$--modules, we have
$B\cong  H^0(G\pi A)$.
By the above we know that
\[
\CD
\Hom_A^{}(Q_i, B)@>>> \Hom_A^{}(P_i, B)
\endCD
\]
is an isomorphism of right $B$--modules. Applying the 
functor $\Hom_B^{}(-,B)$ to it, and recalling that
\[
\Hom_B^{}\big(
\Hom_A^{}(P, B),B\big) =B\otimes_A^{}P,
\]
we deduce that
\[
\CD
B\otimes_A^{}P_i@>>> B\otimes_A^{}Q_i
\endCD
\]
is also an isomorphism.
\eprf

\lem{Ln=0.3}
Any $\s$--inverting ring homomorphism $A\ri C$ factors through
the natural map $A\ri\ct(\pi A,\pi A)\op$.
\elem

\prf
By Lemma~\ref{LTtoD.1} the functor $D(A)\ri D(C)$, taking $X$
to $C\oti X$, factors as
\[
\CD
D(A) @>\pi>> \ct @>T>> D(C).
\endCD
\]
Hence we have ring homomorphisms
\[
\CD
\Hom_{D(A)}^{}(A,A) @>>> \ct(\pi A,\pi A)
@>>> \Hom_{D(C)}^{}(T\pi A,T\pi A).
\endCD
\]
Now $T\pi A=C\otimes A=C$. Taking opposed rings, we have 
\[
\CD
A @>>> \ct(\pi A,\pi A)\op
@>>> C,
\endCD
\]
and the composite is easily seen to be the given map $A\ri C$.
\eprf

\thm{Tn=0.4}
The natural map $A\ri\ct(\pi A,\pi A)\op$ is the initial object
in the category of $\s$--inverting homomorphisms.
\ethm

\prf
Lemma~\ref{Ln=0.2} tells us that the map is $\s$--inverting,
while Lemma~\ref{Ln=0.3} tells us any $\s$--inverting map
factors through it. We need to prove the factorisation unique.
We will prove the uniqueness even as maps of left $A$--modules.

Assume therefore that we are given a $\s$--inverting ring
homomorphism $A\ri C$. Suppose we have a factorisation,
as maps of left $A$--modules,
\[
\CD
A @>>> \ct(\pi A,\pi A)\op
@>>> C.
\endCD
\]
We wish to show it unique. 

By Lemma~\ref{Ln=0.2}, the map $A\ri \ct(\pi A,\pi A)\op$
can be identified with $A\ri H^0(G\pi A)$. By Lemma~\ref{L6.1},
$G\pi A\in\cs^{\leq0}$. This allows us to identify, in the derived
category $\cs=D(A)$, the objects $H^0(G\pi A)$ and ${\{G\pi A\}}^{\geq0}$.
In the derived category $D(A)$, our factorisation of $A\ri C$ becomes
\[
\CD
A @>>> {\{G\pi A\}}^{\geq0}
@>>> C.
\endCD
\]
We can factor this further as
\[
\CD
A @>>> G\pi A @>>> {\{G\pi A\}}^{\geq0}
@>>> C.
\endCD
\]
In the distinguished triangle
\[
\CD
{\{G\pi A\}}^{\leq-1} @>>> G\pi A @>>> {\{G\pi A\}}^{\geq0}
@>>> \T {\{G\pi A\}}^{\leq-1}
\endCD
\]
we have that both ${\{G\pi A\}}^{\leq-1}$ and 
$\T {\{G\pi A\}}^{\leq-1}$ lie in
$\cs^{\leq-1}$. Since $C\in \cs^{\geq0}$, we conclude that the 
map 
\[
\CD
\cs({\{G\pi A\}}^{\geq0},C)
@>>> \cs(G\pi A, C)
\endCD
\]
is an isomorphism. The factorisation 
\[
\CD
A @>>> G\pi A @>>> {\{G\pi A\}}^{\geq0}
@>>> C
\endCD
\]
is completely determined by 
\[
\CD
A @>>> G\pi A 
@>>> C.
\endCD
\]
Now consider the distinguished triangle
\[
\CD
k @>>> A @>>> G\pi A 
@>>> \T k.
\endCD
\]
We know that $k\in\car$. But then, for every $n\in\zz$,
\[
\begin{array}{rcl}
\Hom_{D(A)}^{}(\T^n k,C) &=& \Hom_{D(C)}^{}(\T^n C\oti k\,\,,\,C) \\
            &=& 0.
\end{array}
\]
The last equality is because $C\oti k=0$, by Lemma~\ref{LTtoD.1}.
 From the distinguished triangle we conclude that 
\[
\CD
\cs(G\pi A, C) 
@>>> \cs(A,C)
\endCD
\]
is an isomorphism. The factorisation
\[
\CD
A @>>> G\pi A 
@>>> C
\endCD
\]
is unique.
\eprf

\section{The case $n<0$}
\label{n<0}

In this section we will study what happens to
the groups $\ct(\pi A,\T^n\pi A)=H^n(G\pi A)$ when $n<0$. We will prove that
they vanish if and only if $\rs$ is stably flat; that is,
if and only if $\Tor_n^A(\rs,\rs)=0$ for all $n>0$. We even prove more. We prove that the first non-vanishing $\Tor_n^A(\rs,\rs)$ is isomorphic (up to
changing $n$ to $n-1$) with the first non-vanishing $H^{-n}(G\pi A)$.
For the precise statement see Theorem~\ref{C8.6}.

It might help to give a sketch of the argument. Lemmas~\ref{Ln<0.1}
and \ref{Ln<0.2}
prove that all the homology groups $H^n(G\pi A)$ are naturally modules 
over $\rs$. Lemma~\ref{L8.5} is a technical observation: let $A\ri B$ be a ring
homomorphism, and let $M$ be a $B$ module. Under some hypotheses 
one can say something about
$\Tor^A_i(B,M)$. The idea is
to apply these observations
in the case where $B=\rs$ and the $\rs$--modules in question
are $H^n(G\pi A)$.

Lemma~\ref{Ln<0.3} and Remark~\ref{Rn<0.4} are the crucial part of
the argument. They introduce the spectral sequence which does the
work. Lemma~\ref{Ln<0.5} tells us that the technical conditions
of Lemma~\ref{L8.5} are satisfied in the case of
the ring homomorphism $A\ri\rs$. And then Theorem~\ref{C8.6}
clinches the computation.

\lem{Ln<0.1}
Let $M$ be any $A$--module. There is an isomorphism
of left $A$--modules 
\[
H^0(G\pi M) \cong \br\otimes_A^{} M.
\]
The $A$--module structure on $H^0(G\pi M)$ therefore
extends (uniquely) to an $\rs$--module structure.
\elem

\prf
Put $B=\rs$, and $\theta\co A\ri B$ the initial
$\s$--inverting homomorphism.
Let $x\in D(A)$ be any object. The map
\[
\CD
A\oti x @>\theta\otimes 1>>B  \oti x
\endCD
\]
gives a natural transformation from the identity functor to
$B\oti(-)$. In the distinguished triangle
\[
\CD
k@>>>  x @>\eta_x^{}>> G\pi x @>>> \T k
\endCD
\]
the object $k$ lies in $\car$. Applying the functor
$\Hom_{D(A)}^{}\big(-\,,\, B  \oti x  \big)$ to this
triangle, we deduce an exact sequence
\[
\CD
\Hom_{D(A)}^{}\big(\T k\,,\, B  \oti x  \big) @. \\
@VVV @. \\
\Hom_{D(A)}^{}\big(G\pi x\,,\, B  \oti x  \big)
@>>>
\Hom_{D(A)}^{}\big(x\,,\, B  \oti x  \big) \\
@. @VVV \\
@. \Hom_{D(A)}^{}\big(k\,,\, B  \oti x  \big).
\endCD
\]
But for every $n\in\zz$ there is an isomorphism
\[
\Hom_{D(A)}^{}\big(\T^n k\,,\, B  \oti x  \big)=
\Hom_{D(B)}^{}\big(\T^n B\oti k\,,\, B  \oti x\big).
\]
This vanishes since, by Lemma~\ref{LTtoD.1},
$B\oti k=0$. From the exact sequence
we conclude that the map
\[
\CD
\Hom_{D(A)}^{}\big(G\pi x\,,\, B  \oti x\big)@>>>  
\Hom_{D(A)}^{}\big(x\,,\, B  \oti x\big)
\endCD
\]
is an isomorphism. Our map $\theta\otimes 1\co A\oti x\ri B\oti x$ factors
uniquely as
\[
\CD
x@>\eta_x^{}>> G\pi x @>\p_x^{}>>B  \oti x.
\endCD
\]
The uniqueness allows us to easily show that the $\p_x^{}$ assemble to
 a natural
transformation $\p\co G\pi(-)\Longrightarrow B\oti(-)$. 
Applying the functor $H^0$,
we have a natural transformation
\[
\CD
H^0(\p)\co H^0\big(G\pi(-)\big)\Longrightarrow H^0\big(B\oti(-)\big).
\endCD
\]
What we will show is that,when $x$ is a chain complex concentrated in
degree~0 (ie $x$ is just a module), then $H^0(\p)$ is
an isomorphism. Observe that, when $x$ is just a 
module concentrated in degree 0, then
$H^0\big(B\oti(-)\big)$ simplifies to $B\otimes_A^{}x$.

Let $x$ be $A$, viewed as an object in $D(A)$ concentrated
in degree 0. We have maps
\[
\CD
A@>\eta_A^{}>> G\pi A @>\p_A^{}>>B  \oti A=B.
\endCD
\]
Applying the functor $H^0$, this becomes
\[
\CD
A@>H^0(\eta_A^{})>> H^0(G\pi A) @>H^0(\p_A^{})>>B; 
\endCD
\]
By Lemma~\ref{Ln=0.1} and Theorem~\ref{Tn=0.4} there is an isomorphism
of left $A$--modules $\rho\co H^0(G\pi A)\ri B$, so that the composite
\[
\CD
A@>H^0(\eta_A^{})>> H^0(G\pi A) @>\rho>>B
\endCD
\]
equals $\theta\co A\ri B$. But in the proof of  Theorem~\ref{Tn=0.4}
we saw that any such factorisation (as maps of left $A$--modules) 
is unique. Hence $\rho=H^0(\p_A^{})$,
and $H^0(\p_A^{})$ must be an isomorphism. 
Because both $H^0\big(G\pi(-)\big)$ and 
$B\otimes_A^{}(-)$ commute with direct sums, $H^0(\p_x^{})$ must be an 
isomorphism for any free $A$--module $x$.

Let $M$ be any $A$--module. Choose a free module $F$ surjecting 
onto $M$. We have a short exact sequence of $A$--modules
\[
\CD
0 @>>> K @>>> F @>>> M @>>> 0.
\endCD
\]
We deduce a commutative diagram with exact rows
\[
\CD
H^0(G\pi K) @>>> H^0(G\pi F) @>>> H^0(G\pi M) @>>> H^1(G\pi K) \\
@V{H^0(\p_K^{})}VV  @V{H^0(\p_F^{})}VV @V{H^0(\p_M^{})}VV @VVV \\
B\otimes_A^{}K  @>>> B\otimes_A^{}F @>>> B\otimes_A^{}M @>>> 0.
\endCD
\]
By Lemma~\ref{L6.1} we know that $G\pi K\in D(A)^{\leq0}$, hence
$H^1(G\pi K)=0$. By the above we know that $H^0(\p_F^{})$ is
an isomorphism. This allows us to conclude first that
$H^0(\p_M^{})$ is surjective. This being true for every 
$A$--module $M$, it must be true for $K$. This means
$H^0(\p_K^{})$ is surjective, and hence $H^0(\p_M^{})$ 
must be an isomorphism.

We have proved that $H^0(G\pi M)$ is isomorphic to $B\otimes_A^{}M$,
which is obviously a module over $B=\rs$. The fact that the $B$--module
structure is unique is easy: To say that $X$ is an $A$--module is
to give a ring homomorphism
\[
\CD
A @>>> \Hom_\zz^{}(X,X).
\endCD
\]
To say this extends to a $B$--module structure is to give a 
factorisation of the ring homomorphism through $\theta\co A\ri B$.
The fact that $\theta$ is initial says that any such factorisation
is unique.
\eprf

\lem{Ln<0.2}
For any local object $x\in D(A)$ and any integer $n\in\zz$,
the $A$--module structure  on $H^n(x)$ extends (uniquely)
to a $\rs$--module structure.
\elem

\prf
Replacing $x$ by $\T^n x$, we may assume $n=0$. The object
$x$ is local. Lemma~\ref{Gtrunc} tells us that so
is ${\{x^{\leq0}\}}^{\geq0}$. That is the module $M=H^0(x)$,
viewed as a complex concentrated in degree 0, is a local
object. We need to prove that $M$ is a module over $\rs$.

Consider the distinguished triangle
\[
\CD
k @>\alpha>> M @>>> G\pi M @>>> \T k.
\endCD
\]
Since $M$ is local and $k\in\car$, it follows that $\alpha\co k\ri M$
must vanish. Therefore the map $G\pi M\ri \T k$ admits
a splitting; there is a split inclusion $\beta\co\T k\ri G\pi M$.
By Reminder~\ref{Rj^*} the object $G\pi M$ is local,
while $\T k\in\car$. It follows that  $\beta\co\T k\ri G\pi M$
must vanish. Hence $k=0$, and $M$ is isomorphic to $G\pi M$.
But then
\[
M=H^0(M)= H^0(G\pi M),
\]
and by Lemma~\ref{Ln<0.1} $H^0(G\pi M)$ is naturally a $\rs$--module.
\eprf

The next lemmas are
based on studying two hyperTor spectral sequences. The
key one,
of Remark~\ref{Rn<0.4}, has for its 
$E_2$ term $\Tor_{-i}^A\big(\rs\,,\,H^j(G\pi A)\big)$. Lemma~\ref{Ln<0.2}
tells us that $H^j(G\pi A)$ is naturally a $\rs$--module. Thus we are
interested in general lemmas that apply to $\Tor_{n}^A(B,M)$, where 
$A\ri B$ is a ring homomorphism and $M$ is a $B$--module. We do this by means 
of another hyperTor spectral sequence.

\lem{L8.5}
Suppose $A \ri B$ is a ring homomorphism such that the multiplication map
$\mu\co B\otimes_A^{}B\ri B$ is an isomorphism.  
Suppose also that for some $n \geq 1$
\[
\text{\rm Tor}^A_i(B,B)~=~0~(1\leq i\leq n).
\]
Then for every $B$--module $M$ we have:
\be
\item
The multiplication map $B\otimes_A^{}M\ri M$ is an isomorphism.
\item
$\text{\rm Tor}^A_i(B,M)=0$ for all $1\leq i\leq n$.
\ee
\elem

\prf
Choose a resolution of $M$ by free (left) $B$--modules
\[
\CD
\cdots @>>> Q^{-2} @>>> Q^{-1} @>>> Q^0 @>>> M @>>> 0,
\endCD
\]
and choose a resolution of $B$ by free right $A$--modules
\[
\CD
\cdots @>>> P^{-2} @>>> P^{-1} @>>> P^0 @>>> B @>>> 0.
\endCD
\]
The tensor product $P\otimes_A^{} Q$ gives a double
complex whose cohomology computes $\text{\rm Tor}^A_{-i-j}(B,M)$.
But there is a spectral sequence for it, whose $E_1$ term is
\[
E_1^{i,j}~=~\text{\rm Tor}^A_{-j}(B,Q^i).
\] 
Now $E_1^{i,0}=B\otimes_A^{}Q^i=Q^i$, since $Q^i$ is free
and, by hypothesis,  $B\otimes_A^{}B\ri B$ is an
isomorphism. In $E_2$, we have
\[
E_2^{i,0}~=~\left\{
\begin{array}{ccl}
M &\qquad & \text{if $i=0$}\\
0 &\qquad & \text{otherwise.}
\end{array}
\right.
\]
By hypothesis, we also have $\text{\rm Tor}^A_{-j}(B,B)=0$ for all
$1\leq -j\leq n$. Since $Q^i$ are free, this gives $\text{\rm
Tor}^A_{-j}(B,Q^i)=0$, for all $i$ and for all $1\leq -j\leq n$.  In
other words, $E_1^{i,j}=0$ if $1\leq -j\leq n$, and hence $E_2^{i,j}=0$
if either $j=0$, $i\neq0$, or if $1\leq -j\leq n$.  The assertions of
the lemma immediately
follow.
\eprf

\lem{Ln<0.3}
The map 
\[
\CD
\br\otimes_A^{}A @>1\oti\eta_A^{}>> \br\oti G\pi A
\endCD
\]
is an isomorphism.
\elem

\prf
Consider the distinguished triangle
\[
\CD
k @>>> A @>\eta_A^{}>> G\pi A @>>> \T k.
\endCD
\]
We know that $k\in\car$, and Lemma~\ref{LTtoD.1} tells
us that $\br\oti k=0$. Tensoring the triangle with $\rs$,
the lemma follows.
\eprf

\rmk{Rn<0.4}
Lemma~\ref{Ln<0.3} produced an isomorphism $\rs=\br\oti G\pi A$
in the derived category $D(A)$. Of course, there is a spectral sequence
which computes the cohomology of $\br\oti G\pi A$. The $E_2$ term
is given by
\[
E_2^{i,j}=\Tor_{-i}^A\big(\rs\,,\, H^j(G\pi A)\big).
\]
Lemma~\ref{Ln<0.3} can be viewed as telling us the limit 
of this spectral sequence. 
In the rest of the section we will study the consequences.
\ermk

\lem{Ln<0.5} We have:
\be
\item
The multiplication map $\mu\co\br\otimes_A^{}\br\ri\br$ is an isomorphism.
\item $\Tor_1^{A}(\rs,\rs)=0$.
\ee
\elem

\prf
The results of Lemma~\ref{Ln<0.5} are not new. They first appeared in an
article by Bergman and Dicks~\cite{Bergman-Dicks78}. The proof that
$\mu$ is an isomorphism may be found in 
\cite[(4) on page 298]{Bergman-Dicks78}, combined with the remark in 
Construction~2.2 on page~300. The vanishing of $\Tor_1^{A}(\rs,\rs)$
is in \cite[(95) on page 326]{Bergman-Dicks78}. See also
Schofield~\cite[page 58]{Schofield}. Even though the results
are known, both statements are
essentially immediate from the spectral sequence of
Remark~\ref{Rn<0.4}. Hence we include the proof.

By Lemma~\ref{L6.1} we know that $G\pi A\in D(A)^{\leq 0}$. This means
$H^j(G\pi A)=0$ when $j>0$. For $\Tor_{-i}^A(-,-)$, we know it vanishes
whenever $i>0$. In the spectral sequence of Remark~\ref{Rn<0.4}
we therefore have $E_2^{i,j}=0$ unless $i\leq0$ and $j\leq0$. The
spectral sequence is third quadrant. 

This immediately means that $E_2^{0,0}=E_\infty^{0,0}$ and
$E_2^{-1,0}=E_\infty^{-1,0}$. 
Lemma~\ref{Ln<0.3} tells us precisely that the map 
\[
\CD
1\oti\eta_A^{}\co\br\otimes_A^{}A @>>> \br\oti G\pi A
\endCD
\]
is an isomorphism. Evaluating $H^{-1}$ of this isomorphism,
we have that $E_\infty^{-1,0}$ is a quotient of $H^{-1}(\rs)=0$,
and hence
\[
E_2^{-1,0}=\Tor_1^{A}(\rs,\rs)=0.
\]
Evaluating $H^0$ of the isomorphism of Lemma~\ref{Ln<0.3}
and recalling that $E_2^{0,0}=E_\infty^{0,0}$, we have
an isomorphism
\[
\CD
1\otimes_A^{} H^0(\eta_A^{})\co\br\otimes_A^{}A 
@>>> \br\otimes_A^{} H^0(G\pi A).
\endCD
\]
By Lemma~\ref{Ln=0.1} we know that the map $H^0(\eta_A^{})\co A\ri
H^0(G\pi A)$ can be identified with the natural homomorphism
$\theta\co A\ri\rs$. It follows that
\[
\CD
1\otimes_A^{} \theta\co\br\otimes_A^{}A @>>> \br\otimes_A^{} \br
\endCD
\]
is an isomorphism. But the composite
\[
\CD
\br\otimes_A^{}A @>1\otimes_A^{} \theta>> \br\otimes_A^{} \br
@>\mu>> \rs
\endCD
\]
is clearly the identity. This makes $\mu$ left-inverse to the 
invertible map $1\otimes_A^{} \theta$. Hence $\mu$ must be the
two-sided inverse, and is invertible.
\eprf

\thm{C8.6}
Suppose $\text{\rm Tor}^A_{i}(\rs,\rs)=0$, for all $1\leq i\leq n$.
Then for all $1\leq i\leq n-1$ we have $H^{-i}(G\pi A)=0$,
and
\[
\text{\rm Tor}^A_{n+1}(\rs,\rs)~=~H^{-n}(G\pi A).
\]
\ethm

\proof
The proof is a slightly more sophisticated computation 
with the spectral
sequence of Remark~\ref{Rn<0.4}. Recall that we have
a spectral sequence whose $E_2$ term is
\[
E_2^{i,j}~=~\text{\rm Tor}^A_{-i}\big(\rs,H^j(G\pi A)\big),
\]
which converges to $H^{i+j}(\rs)$. By Reminder~\ref{Rj^*} 
the object $G\pi A$
is local. By Lemma~\ref{Ln<0.2} the homology
groups $H^j(G\pi A)$ are all modules over $\rs$.
By Lemma~\ref{Ln<0.5} we know that the multiplication
map $\mu\co\br\otimes_A^{}\br\ri\rs$ is an isomorphism.
Lemma~\ref{L8.5} now applies, and we deduce that
if $1\leq -i\leq n$ then $E_2^{i,j}=0$. This forces
the differential
\[
E_2^{-i-1,0}\ri E_2^{0,-i}
\]
to be an isomorphism, for all $1\leq i\leq n$. For 
$1\leq i\leq n-1$ we read off that $H^{-i}(G\pi A)=0$.
For $i=n$, we deduce that
\[
\text{\rm Tor}^A_{n+1}(\rs,\rs)~=~H^{-n}(G\pi A). \eqno{\Box}
\]

\end{document}